\documentclass[twoside, 12pt]{amsart}
\usepackage[active]{srcltx} 
\usepackage[all]{xy}
\CompileMatrices
\usepackage{ifthen}
\usepackage{amsfonts}
\usepackage{amssymb}
\usepackage{amsthm}
 \newtheorem{thm}{Theorem}[section]
 \newtheorem{prop}[thm]{Proposition}
 \newtheorem{cor}[thm]{Corollary}
 \newtheorem{lem}[thm]{Lemma}
  \newtheorem{ques}[thm]{Question}
 
 \newtheorem*{bthm}{Theorem}

\theoremstyle{definition}
\newtheorem{defn}[thm]{Definition}

\theoremstyle{remark}
\newtheorem{rem}[thm]{Remark}

  1

\newcommand{\fp}{\ifmmode {\mathbb{F}_p}\else$\mathbb{F}_p$\ \fi}
\newcommand{\zp}{\ifmmode {\mathbb{Z}_p}\else$\mathbb{Z}_p$\ \fi}
\newcommand{\z}{\mathbb{Z}}
\newcommand{\zpMod}{\ifmmode\mbox{$\zp$-Mod}\else$\zp$-Mod \fi}
\newcommand{\Mod}{\ifmmode\mbox{$\Lambda$-Mod}\else$\Lambda$-Mod \fi}
\renewcommand{\mod}{\ifmmode\mbox{$\Lambda$-mod}\else$\Lambda$-mod
\fi}
\newcommand{\La}{\ifmmode{\Lambda}\else{$\Lambda$}\fi}
\newcommand{\rhom}{\mbox{$\mathbf{R}\mbox{Hom}$}}
\newcommand{\Hom}{{\mathrm{Hom}}}
\newcommand{\Ext}{{\mathrm{Ext}}}
\newcommand{\Tor}{{\mathrm{Tor}}}
\newcommand{\tor}{{\mathrm{tor}}}
\newcommand{\E}{{\mathrm{E}}}

\newcommand{\pd}{{\mathrm{pd}}}
\newcommand{\cd}{{\mathrm{cd}}}
\newcommand{\rk}{{\mathrm{rk}}}
\newcommand{\depth}{{\mathrm{depth}}}
\newcommand{\K}{\mathbf{K}}
\newcommand{\D}{\mathbf{D}}
\newcommand{\Lotimes}{\otimes^\mathbf{L}}
\renewcommand{\H}{\mathrm{H}}
\newcommand{\du}{\mathrm{D}}
\renewcommand{\t}[1]{\mbox{$\mbox{\rm T}_{#1}$}}
\newcommand{\fa}{\mbox{\ for all \ }}
\newcommand{\M}{\ifmmode {\frak M}\else${\frak M}$ \fi}
\newcommand{\m}{\ifmmode {\frak m}\else$\frak m$ \fi}
\newcommand{\p}{\ifmmode {\frak p}\else${\frak p}$\ \fi}
\renewcommand{\P}{\ifmmode {\frak P}\else${\frak P}$\ \fi}
\newcommand{\e}{\ifmmode {\mathcal{E}}\else$\mathcal{E}$ \fi}
\newcommand{\n}{\mbox{${\Bbb N}$}}

\newcommand{\C}{\mbox{$\mathcal{ C}$}}

\newcommand{\G}{\ifmmode {\mathcal{G}}\else${\mathcal{G}}$\ \fi}
\renewcommand{\d}{\ifmmode {\mathcal{ D}}\else${\mathcal{D}}$\ \fi}
\newcommand{\A}{\ifmmode {\mathcal{A}}\else${\mathcal{ A}}$\ \fi}

\renewcommand{\projlim}[1] {{\lim\limits_{\stackrel{\displaystyle
\longleftarrow}{#1}}}}
\newcommand{\projlimsc}[1]
{{\lim\limits_{\stackrel{\scriptstyle \longleftarrow}{#1}}}}
\newcommand{\dirlim}[1]
{{\lim\limits_{\stackrel{\displaystyle \longrightarrow}{#1}}}}
\renewcommand{\in}{\ \epsilon\ }

\newcommand{\kl}{[\![}
\newcommand{\kr}{]\!]}
\newcommand{\Qp}{\ifmmode {{\Bbb Q}_p}\else${\Bbb Q}_p$\ \fi}
\newcommand{\qp}{\ifmmode {{\Bbb Q}_p}\else${\Bbb Q}_p$\ \fi}
\newcommand{\Q}{\ifmmode {\Bbb Q}\else${\Bbb Q}$\ \fi}

\newcommand{\Ind}{\mathrm{Ind}}
\newcommand{\ho}{\mathrm{Ho}}

\begin{document}

\title[On the structure theory of the Iwasawa algebra]{On the structure theory of the Iwasawa algebra of a $p$-adic Lie group}%
\author{Otmar Venjakob}%
\address{Universit\"{a}t Heidelberg\\ Mathematisches Institut\\
 Im Neuenheimer Feld 288\\
 69120 Heidelberg, Germany.}
\email{otmar@mathi.uni-heidelberg.de}
\urladdr{http://www.mathi.uni-heidelberg.de/\textasciitilde
otmar/}



\begin{abstract}
This paper is motivated by the question whether there is a nice
structure theory of finitely generated modules over the  Iwasawa
algebra, i.e.\ the completed group algebra, \La\ of a $p$-adic
analytic group $G.$ For $G$ without any  $p$-torsion element we
prove that \La\ is an Auslander regular ring. This result enables
us to give a good definition of the notion of a  {\em pseudo-null}
\La-module. This is classical when $G=\mathbb{Z}^k_p$ for some
integer $k\geq 1,$ but was previously unknown in the
non-commutative case. Then the category of \La-modules up to
pseudo-isomorphisms is studied and we obtain a weak structure
theorem for the $\zp$-torsion part of a finitely generated
\La-module. We also prove a local duality theorem and a version of
Auslander-Buchsbaum equality. The arithmetic applications to the
Iwasawa theory of abelian varieties are published elsewhere.
\end{abstract}

\maketitle

\PUSH{intro.tex}%
\section*{Introduction}

The increasing interest in the Iwasawa algebra, i.e.\ the
completed group algebra
 \[\Lambda(G)=\zp\kl G\kr\]
 of a compact $p$-adic analytic group $G$
stems at least from the following  two reasons: Firstly, in a
purely algebraic sense $\La(G)$ is a prototype of a ``regular"
 - in general non-commutative - ring. It seems very likely that the category
  of finitely generated $\La(G)$-modules, which we denote by \mbox{$\La(G)$-mod,} has a very rich structure,
  and our aim is to give some initial evidence in this direction.  Secondly,
 compact $p$-adic Lie groups occur  naturally in arithmetic geometry as the Galois groups
 of certain infinite  extensions $k_\infty$ of a number field $k,$
 for example, $k_\infty$ could be taken to be the field obtained
by adjoining to $k$ the coordinates of all $p$-power division
points on some abelian variety defined over
 $k.$ It is then natural to study certain  arithmetic objects
 defined over $k_\infty$ (for example, the $p$-Hilbert class field
 of $k_\infty$ or the Pontryagin dual of the Selmer  group of an
 abelian variety over $k_\infty$)  as modules over the Iwasawa algebra $\La(G)$  (cf.\ \cite{harris},
 \cite{coates-howsonII}, \cite{coates}).
Though this last aspect widely motivated our study of $\La(G)$ the
results of this paper are ``purely algebraic". For a detailed
discussion of its arithmetic applications we refer the reader to
\cite{ochi-ven} and a forthcoming paper \cite{ven}.\\ Let us
assume for simplicity that the $p$-adic analytic group $G$ is
torsion-free and a pro-$p$-group. Then, $\La(G)$  is a (both left
and right) Noetherian ring without zero-divisors. Furthermore, by
results of Brumer it is known that $\La(G)$ has finite projective
dimension equal to $\pd (\La(G))= \dim (G) +1,$ where $\dim(G)$
denotes the dimension of $G$ as $p$-adic analytic manifold and
agrees with its $p$-cohomological dimension. Thus in some sense
$\La(G)$ can be considered as a ``regular" ring, and it is natural
to ask if there could exist  a  structure theory for finitely
generated $\La(G)$-modules, which is parallel to the classical
theory (studied in detail in Bourbaki) when $G=\mathbb{Z}_p^k$ for
some integer $k\geq1.$ A first basic step in this direction is to
answer the following question
 \begin{quote}
 What is a good definition of  {\em pseudo-null} \/ resp.\ {\em
 pseudo-isomorphism} in the context of $\La(G)$-modules?
  \end{quote}
We recall that for  a commutative Noetherian ring $R$ and  a
finitely generated $R$-module $M$ the definition is the following:
The dimension of $M$ is defined to be the Krull dimension of the
support of $M$ in $\mathrm{Spec}(R)$ and $M$ is said to be
pseudo-null, if its codimension is greater than 1.  M.\ Harris
\cite[1.12]{harris} already proposed a vague definition of
pseudo-null using a certain filtration of $\La(G),$ which in
general differs from the $\M$-adic one, where $\M$ denotes the
maximal ideal of $\La(G),$ and cannot be described easily. Besides
some more or less trivial cases it turned out very difficult to
verify whether a module is
pseudo-null.\\

 In this article we  give an  answer to this
question using the following strategy. First, we establish a
dimension theory for finitely generated $\La(G)$-modules. Our
approach here has been inspired by the work  of Bj\"{o}rk \cite{erik},
who showed that each finitely generated module over an {\em
Auslander regular} or more generally {\em Auslander Gorenstein}
ring (for the definitions see \ref{defAuslander}) is intrinsically
equipped with a canonical filtration
 $$ \t0 (M) \subseteq  \t1 (M)\subseteq\cdots\subseteq\t{d-1}
(M)\subseteq\t{d} (M)=M.$$ Using this  filtration he defines the
dimension of a \La-module $M.$  It turns  out that for a
commutative regular local ring this dimension equals the Krull
dimension and that $\t{i}(M)$ is just the maximal submodule of $M$
with dimension less or equal to $i.$\\
 Thus the following theorem states a fundamental structure property of $\La(G)$  and will be
crucial for the applications in Iwasawa theory ( \cite{ven}). Here
and in what follows we assume that $G$ is a compact $p$-adic
analytic group with no element of order $p$ (but not necessarily
pro-$p$).

\begin{bthm} (Theorem \ref{Lambda-Auslanderreg})
$\La(G)$ is an Auslander regular ring.
\end{bthm}

For the purpose of studying the $\zp$-torsion part of
$\La(G)$-module the following consequence for the completed group
algebra $\mathbb{Z}/p^n\kl  G\kr \cong\La/p^n$ with coefficients
in $\mathbb{Z}/p^n$  becomes very useful.

\begin{bthm} (Theorem \ref{fp-Auslander})
 \begin{enumerate}
\item $\zp/p^m\kl  G\kr $ is an Auslander-Gorenstein ring with injective
dimension equal to $\cd_p(G).$
\item $\fp\kl  G\kr $ is an Auslander regular ring of  dimension
$\cd_p(G).$
\end{enumerate}
\end{bthm}

If one is ready to neglect the $\zp$-torsion part - which can be
controlled by the last theorem and a further result below, anyway
- then it might be convenient to consider also modules over the
ring $\qp\kl G\kr:=\zp\kl G\kr\otimes_\zp\qp.$

\begin{bthm} (Theorem \ref{qp-Auslander})
 $\qp\kl G \kr$ is an Auslander regular ring
of  dimension $\cd_p(G).$
\end{bthm}

Using these results, it is quite obvious how to define
pseudo-null:
 \begin{quote}
 A finitely generated \La-module is called {\em pseudo-null} if
 and  only if its co-dimension is  greater or equal to $2.$
 \end{quote}
In the case $G=\mathbb{Z}_p^k$ this is just the usual definition.
So we are convinced that our definition is the right
generalization to the non-commutative case.

Once having available the notation of pseudo-null modules it is
natural to ask whether there holds a structure theorem of
$\La(G)$-modules {\em up to pseudo-\-iso\-morphism} which is
analogous to the well known structure theorem for finitely
generated modules over a commutative regular ring.

The first result in this direction was obtained by the author in
his thesis \cite{ven-diss}, and determines the $\La(G)$-structure
of the $\zp$-torsion submodule of a finitely generated
$\La(G)$-module  in the quotient category of $\mod$ by the Serre
subcategory  of pseudo-null \La-modules. We write $\mathcal{PN}$
for the full subcategory of $\mod$ consisting of the pseudo-null
modules.

\begin{bthm} (Theorem \ref{p-structure})
Assume that $G$ is  a $p$-adic analytic group without $p$-torsion
such that both $\La=\La(G)$ and $\La/p$ are integral.  Then, for
any \La-module $M,$ there exist uniquely (up to order) determined
natural numbers $n_1,\ldots,n_r$ such that
 $$\tor_\zp M\equiv \bigoplus_{1\leq i\leq r} \La/p^{n_i}\ \mbox{\rm  mod
 }{\mathcal{
 PN}}.$$
\end{bthm}

Both the hypotheses that $\La(G)$ and $\La(G)/p\La(G)$ have no
zero divisors are known to be true if $G$ is a uniform, pro-$p$
$p$-adic Lie group with no element of order $p$ (\cite{dsms}). Our
results have inspired various  mathematicians to generalize this
result and  search for a full analogue of the structure theory.
Firstly, Susan Howson has generalized it to the submodule of $M$
annihilated by some power of a fixed prime element in the centre
of $\La(G),$ obtaining the stronger result that a
pseudo-isomorphism exists in the category $\mod$
(\cite{howson2001}). Secondly, J. Coates and R. Sujatha
(\cite{co-su1}) and  P. Schneider \cite{schneider2001} have proven
that, in the quotient category $\mod/\mathcal{
 PN},$ every finitely generated
torsion $\La(G)$-module is isomorphic to a  direct sum of cyclic
$\La(G)$-modules, when $G$ is any extra-powerful pro-$p$ $p$-adic
Lie group with no element of order $p.$ Schneider has shown that
this theorem is in fact a consequence of a general result of
Chamarie \cite{cham} on the structure of modules over
non-commutative Krull domains, whereas Coates and Sujatha show
that it can also be proven using techniques from the algebraic
theory of micro-localization. It is still unknown whether  the
left cyclic modules which occur can be chosen to be quotients of
$\La(G)$ by left principal ideals.

Now we want to state two further main results on the structure of
$\La(G),$ if $G$ is a pro-p Poincar\'{e} group  of finite
cohomological dimension and such that $\La=\La(G)$ is Noetherian.
The first result tells us that $\La(G)$ ``admits local duality \`{a}
la Grothen\-dieck", i.e.\ if local cohomology is defined in an
natural way (see section \ref{localdualitysec}), we obtain

\begin{bthm} (Theorem \ref{localduality})
For any $M\in\mbox{$\La(G)$-mod},$ there are canonical
isomorphisms
 $$\E^i(M)\cong\Hom_\La(\H^{d-i}_\M(M),\H^{d}_\M(\La))\cong\H^{d-i}_\M(M)^\vee,$$
where $d=\cd_p(G)+1.$
\end{bthm}

The second result generalizes the Auslander-Buchsbaum equality.

\begin{bthm} (Theorem \ref{ab-equa})
For any $M\in\mod$, it holds
 $$\pd_\La(M)+\depth_\La(M)=\depth_\La (\La).$$
\end{bthm}

\vspace{0.5cm}

 \textsc{Acknowledgements.}\\

This article is based on  a part of my Dissertation, Heidelberg
2000, and I would like to thank my supervisor Kay Wingberg most
warmly for leading me to the nice field of ``higher dimensional"
Iwasawa theory. Without his advice and confidence this project
would not have been possible.  I am also much indebted to John
Coates and Alexander Schmidt, whose comments on my manuscript
helped improve the exposition. Peter Schneider is heartily
acknowledged for pointing out some inaccuracies.

%
\POP
\PUSH{basics.tex}%

\section{Basic properties of \La-modules}
\subsection{Preliminaries}

 The aim  of this work
is to give some complements to the theory of $\Lambda$-modules,
where we denote by $\Lambda=\Lambda(G)$ the completed group
algebra of a profinite group $G$ over $\zp$
 $$\Lambda(G)=\zp\kl G\kr =\projlim{U}\zp[G/U].$$
Here $U$ runs through the open normal subgroups of $G$. We start
by recalling some well-known  facts concerning $\Lambda$,
 proofs
of which can be found  in \cite[V\S2]{nsw}. By a (left)
$\La$-module $M$ we understand a separated topological module, i.e
$M$ is a Hausdorff topological $\zp$-module with a continuous
$G$-action. Since the involution of $\La$ given by passing to the
inverses of group elements induces a natural equivalence between
the categories of left and right $\La$-modules, we will often
ignore the difference without further mention. The category
$\C=\C(G)$ of compact $\La$-modules and the category $\d=\d(G)$ of
discrete $\La$-modules will be of particular interest. Both are
abelian categories, and Pontryagin duality defines a contravariant
equivalence of categories between them. Hence, while $\C$ has
sufficiently many projectives and exact inverse limits the dual
statement holds for $\d.$ \\ By $I_G$ we denote the augmentation
ideal of $\La,$ i.e.\ the kernel of the canonical epimorphism
 $$\zp\kl G\kr \twoheadrightarrow\zp$$
and by $$M_G=M/I_GM$$ the module of coinvariants of $M.$ Then, the
$G$-homology $\H_\bullet(G,M)$ of a compact \La-module $M$ can be
defined as left derived functor of $-_G$ or alternatively as
$\Tor^\La_\bullet(\zp,M),$ where $\Tor$ denotes the left derived
functor of the complete tensor product $-\widehat{\otimes}_\La-.$
We obtain  a canonical isomorphism
$\H_i(G,M)\cong\H^i(G,M^\vee)^\vee$, where $\H^\bullet(G,-)$
denotes the usual $G$-cohomology for a discrete \La-module
considered as a discrete abelian group and $^\vee$ is the
Pontryagin dual.

In order to state the topological Nakayama lemma we define the
radical $\mathrm{Rad}_G$ of \La\ to be the intersection of all
open left maximal ideals. It is a closed two-sided ideal and its
powers define a topology on \La\ which is called the $R$-topology.
If a $p$-Sylow group $G_p$ is of finite index in $G$, then this
topology coincides with the canonical one \cite[5.2.16]{nsw},
$\mathrm{Rad}_G$ is an open ideal of \La\ and all (left) maximal
ideals are open. Furthermore, $\La(G)$ is a local ring if and only
if $G$ is a pro-$p$-group. In this case the maximal ideal of \La\
is equal to $p\La+I_G.$

\begin{lem}(Topological Nakayama Lemma)
\begin{enumerate}
\item If $M\in\C$ and $\mathrm{Rad}_GM=M,$ then $M=0.$
\item Assume that $G$ is a pro-$p$-group, i.e.\ \La\ a local ring
with maximal ideal $\M.$ Then the following holds:
\begin{enumerate}
\item  $M$ is generated
by $x_1,\ldots,x_r$ if and only if $x_i+\M M,$ $i=1,\ldots, r,$
generate $M/\M M$ as $\mathbb{F}_p$-vector space.
\item Let $P\in\C$ be finitely generated. Then $P$ is \La-free if
and only if $P$ is \La-projective.
\end{enumerate}
\end{enumerate}
\end{lem}

Compare \cite{bal-how} for a thorough discussion.

 Concerning the projective dimension $\pd_\La
M$, respectively global dimension $\pd( \La)$ of $\La$, which are
both defined with respect to the category $\C,$ there are the
following results due to Brumer \cite{brumer}, where $\cd_p(G)$
denotes the $p$-cohomological dimension of $G:$
\[\pd_\La \zp= \cd_p(G)\mbox{ and }
 \pd (\La)=\cd_p(G)+1.\]
If $\La$ is Noetherian (e.g. if $G$ is a $p$-adic Lie group, see
below), the forgetful functor from the category $\C$ of compact
\La-modules to the category $\Mod$ of abstract \La-modules defines
an equivalence between the full subcategory $\C_{fg}$ of finitely
generated compact \La-modules and the full subcategory $\mod$ of
finitely generated abstract \La-modules. In particular, the
different definitions of the projective, respectively global
dimension coincide in this case.

\subsection{$p$-adic Lie groups}

In this subsection we recall some facts about (compact) $p$-adic
Lie groups or also called (compact) $p$-adic analytic groups,
i.e.\ the group objects in the category of $p$-adic analytic
manifolds over $\qp.$ There is a famous characterization of them
due to Lazard \cite{la} (see also \cite{dsms} 9.36):

{\em A topological group $G$  is a compact $p$-adic Lie group if
and only if
 $G$ contains a normal open uniformly powerful
pro-$p$-subgroup of finite index.}

Let us briefly recall the definitions:  A pro-$p$-group $G$ is
called {\em powerful}, if  $[G,G]\subseteq G^{p}$ for odd $p,$
respectively $[G,G]\subseteq G^{4}$ for $p=2$ holds. A
(topologically) finitely generated powerful pro-$p$-group $G$ is
{\em uniform} if the $p$-power map induces isomorphisms
 $$P_i(G)/P_{i+1}(G)\stackrel{\cdot^p}{\to}P_{i+1}(G)/P_{i+2}(G),
 \quad i\geq 1,$$
where $P_i(G)$ denotes  the lower central $p$-series  given by
$$P_1(G)=G,\ P_{i+1}(G)=P_i (G)^p[P_i(G),G],\ i\geq 1,$$ (for
finitely generated pro-$p$-groups). It can be shown that for a
uniform group $G$ the sets $G^{p^i}:=\{g^{p^i}|\ g\in G\}$ form
subgroups and in fact $G^{p^i}=P_{i+1}(G),\ i\geq 0.$ For example,
all the congruence kernels of $GL_n(\zp),SL_n(\zp)$ or
$PGL_n(\zp)$ are uniform pro-$p$-groups for $p\neq 2,$ in
particular the lower central $p$-series  of the first congruence
kernel corresponds precisely to the higher congruence kernels. We
should mention also the  following basic result (see \cite{dsms},
p. 62):

{\em A pro-$p$ powerful group is uniform if and only if it has no
element of order $p.$}

It is a remarkable fact that the analytic structure of a $p$-adic
Lie group is already determined by its topological group
structure. Also, the category of $p$-adic analytic groups is
stable under the formation of closed subgroups, quotients and
group extensions (See \cite{dsms}, chapter 10, for these facts).
The following cohomological property is indispensable ( for the
definition of Poincar\'{e} groups see \cite{nsw}).

{\em A $p$-adic Lie group of dimension $d$ (as $p$-adic analytic
manifold) without $p$-torsion is a Poincar\'{e} group at $p$ of
dimension $d.$}

With respect to the completed group algebra we know that $\La(G)$
is Noetherian for any compact $p$-adic Lie group (see \cite{la}V
2.2.4). If, in addition, $G$ is uniform, then $\La(G)$ is an
integral domain, i.e.\ the only zero-divisor in $\La(G)$ is $0$
(\cite{la}). In fact the latter property also holds for any
$p$-adic analytic group without elements of finite order (see
\cite{Ne}). For instance, for $p\geq n+2,$ the group $Gl_n(\zp)$
has no elements of order $p,$ in particular, $GL_2(\zp)$ contains
no elements of finite $p$-power order if $p\geq 5$ (see
\cite{howson} 4.7) .

 In this case (i.e.\ \La\ is both left and
right Noetherian and without zero-divisors) we can form a skew
field $Q(G)$ of fractions of \La\ (see \cite{GoWa}). This allows
us to define the rank of a \La-module:

\begin{defn}
The rank $\rk_\Lambda M$ is defined to be the dimension of
\linebreak $M\otimes_\La Q(G)$ as a left vector space over $Q(G)$
 $$\rk_\Lambda M=\dim_{Q(G)}(M\otimes_\La Q(G)).$$

\end{defn}

\section{Homotopy theory of modules}

In this section we briefly recall some definitions and results
from the homotopy theory of modules  in the setting of U. Jannsen
\cite{ja-is} , who developed further the homotopy theory which was
introduced by Eckmann and Hilton and extended by Auslander and
Bridger \cite{aus}. The proofs can be found in \cite[\S1]{ja-is}
or in \cite[V\S4]{nsw}. Though this theory works in much larger
generality, we restrict ourselves to the case where \La\ is a
Noetherian (= right and left Noetherian)  ring (= associative
unital ring). Furthermore, {\em all \La-modules considered are
assumed to be finitely generated.}

\begin{defn}
A morphism $f:M\to N$ of \La-modules is homotopic to zero, if it
factors through a projective module $P:$
 $$f:M\to P\to N.$$
 Two morphisms $f,g$ are homotopic
$(f\simeq g),$ if $f-g$ is homotopic to zero. Let
$[M,N]=\Hom_\La(M,N)/\{f\simeq 0\}$ be the group of homotopy
classes of morphisms from $M$ to $N,$ and let $\ho(\La)$ be the
category, whose objects are the objects of $\mod$ and whose
morphism sets are given by $\Hom_{\ho(\La)}(M,N)=[M,N],$ i.e.\ the
category of ``\La-modules up to homotopy."
\end{defn}

\begin{rem}
The additive homotopy category of modules is not abelian in
general. It can be shown that it forms a closed model category
(for suitable definitions of (co)fibrations and weak
equivalences). In general, it cannot be extended to a triangulated
category: If it were a triangulated category in general there
would have to exist for any module $M$ a weak equivalence between
$M$ and $\Omega M,$ where $\Omega$ denotes the loop space functor
which will be introduced below. But for a ring \La\ with finite
projective dimension this would imply that all modules in
$\La$-mod are projective.\\ However, if \La\ is a quasi-Frobenius
ring (for the definition and properties see \cite[4.2]{weibel}),
e.g.\ the group algebra of a finite group over a field $\La=k[G],$
 then its
associated homotopy category is triangulated (\cite[IV Ex.
4-8]{gelfand-manin}).
\end{rem}

It turns out that $M$ and $N$ are homotopy equivalent, $M\simeq
N,$ i.e isomorphic in $\ho(\La),$ if and only if $M\oplus P\cong
N\oplus Q$ with projective \La-modules $P$ and $Q.$ In particular,
$M\simeq 0$ if and only if $M$ is projective.\\

\begin{defn}
For $M\in\mod$ we define the higher adjoints of $M$ to be
 $$\E^i(M):=\Ext^i_\La(M,\La),\quad i\geq 0,$$
which are a priori right \La-modules by functoriality and the
right \La-structure of the bi-module \La. By convention we set
$\E^i(M)=0$ for $i<0.$ The \La-dual $\E^0(M)$ will also be denoted
by $M^+.$ If $\La=\La(G)$ is a completed group algebra these
adjoints will be considered as left modules via the involution of
\La.
\end{defn}

\begin{rem}
In case $\La$ is the completed group algebra of a profinite group
$G$ we call the $E^i(M)$ also {\em Iwasawa adjoints} of $M$ as
$E^1(M)$ is isomorphic (up to the inversion of the group action)
to Iwasawa's adjoint $\alpha(M)$ for $G\cong\zp.$
\end{rem}

It can be shown that for $i\geq 1$ the functor $\E^i$ factors
through $\ho(\La)$ defining  a functor
 $$\E^i:\ho(\La)\to\mathrm{mod-}\La.$$

Now we will describe the construction of a contravariant duality
functor, the transpose
 $$\du :\ho(\La)\to\ho(\La),$$
where we denote the homotopy category of right \La-modules also by
$\ho(\La).$
 For every $M\in\mod$  choose a presentation $P_1\to P_0\to
M\to 0$ of $M$ by projectives and define the transpose $\du M$ by
the exactness of the sequence
 $$\xymatrix@1{
    {\ 0 \ar[r] } &  {\ M^+ \ar[r] } &  {\ P_0^+ \ar[r] } &  {\ P_1^+ \ar[r] } &  {\
    \du M\ar[r] }& 0. \
 }$$
Then it can be shown that the functor $\du$ is well-defined and
one has $\du^2=\mathrm{Id}.$ Furthermore, if $\pd_\La M\leq 1$
then $\du M\simeq \E^1(M).$

Sometimes it is also convenient to have the notation of the
$1^{st}$ syzygy or loop space functor $\Omega:\mod\to\ho(\La)$
which is defined as follows (see \cite[1.5]{ja-is}): Choose a
surjection $P\to M$ with $P$ projective. Then $\Omega M$ is
defined by the exact sequence
 $$\xymatrix@1{
    {\ 0 \ar[r] } &  {\ \Omega M \ar[r] } &  {\ P \ar[r] } &  {\ M \ar[r] } &  {\
    0 }. \
 }$$

The next result will be of particular importance:

\begin{prop}\label{canonicalsequ}
For $M\in\mod$ there is a canonical exact sequence
 $$\xymatrix@1{
    {\ 0 \ar[r] } &  {\ \E^1\du M \ar[r] } &  {\ M \ar[r]^{\phi_M }} &  {\ M^{++} \ar[r] } &  {\
    \E^2\du M\ar[r] }& 0, \
 }$$
where $\phi_M $ is the canonical map from $M$ to its bi-dual. In
the following we will refer to the sequence as ``the" canonical
sequence (of homotopy theory).
\end{prop}

A \La-module $M$ is called {\em reflexive} if $\phi_M $ is an
isomorphism from $M$ to its bi-dual $M\cong M^{++}.$

As Auslander and Bridger suggest the module $\E^1\du M$ should be
considered as torsion submodule of $M.$ Indeed, if \La\ is a
Noetherian integral domain this submodule coincides exactly with
the set  of torsion elements $\tor_\La M.$

\begin{defn}\label{tordef}
A \La-module $M$ is called  {\em \La-torsion module} if
$\phi_M\equiv0,$ i.e.\ if $\tor_\La M:=\E^1\du M=M.$ We say that
$M$ is {\em \La-torsion-free} if $\E^1\du M=0.$
\end{defn}

Since $M^{++}$ embeds into a free \La-module (just take the dual
of an arbitrary surjection $\La^m\twoheadrightarrow M^+$) the
torsion-free \La-modules are exactly the submodules of free
modules. A different characterization of torsion(-free) modules
will be given later using dimension theory, see \ref{torchar}.\\

For $\La:=\La(G),$ where  $G$ is a $p$-adic Lie group, the above
definition can be interpreted as follows: {\em A finitely
generated \La-module $M$ is  a \La-torsion module if and only if
$M$ is  a $\La(G')$-torsion module (in the strict sense) for some
open pro-$p$ subgroup $G'\subseteq G$ such that $\La(G')$ is
integral.} Indeed, for any open subgroup $H$ of a $p$-adic Lie
group $G$ there is an isomorphism
$\E^1_{\La(G)}\du_{\La(G)}\cong\E^1_{\La(H)}\du_{\La(H)}$ of
$\La(H)$-modules by  proposition \ref{inducedExt} (ii) (see below)
 and the analogue statement for $\du M$.\\

For a closed subgroup $H\subseteq G$ we denote by
$\Ind^H_G(M)=M\widehat{\otimes}_{\La(H)}\La(G)$ the compact
induction of a $\La(H)$-module to a $\La(G)$-module.

\begin{prop} Let $H$ be a closed\label{inducedExt}
subgroup  of $G$.
\begin{enumerate}
\item For any $M\in\La(H)$-mod and any $i$  we have an
isomorphism of \La-modules
 $$\E^i_{\La(G)}(\Ind_G^H M)\cong \Ind_G^H \E^i_{\La(H)}(M).$$
\item If, in addition, $H$ is an open subgroup, then there is an
isomorphism of $\La(H)$-modules
 $$\E^i_{\La(G)}(M)\cong\E^i_{\La(H)}(M).$$
\end{enumerate}
\end{prop}
\begin{proof}
The first statement is proved in \cite[lemma 5.5]{ochi-ven} while
the second one can be found in \cite[lemma 2.3]{ja-is}.
\end{proof}
\POP
\PUSH{filtrations.tex}%
\section{Auslander regularity}
\subsection{Filtrations of \La-modules}

Since the completed group ring \La\ of a $p$-adic Lie group
without $p$-torsion is (both left and right) Noetherian and has
finite global homological (and therefore finite injective)
dimension we can apply the results of J.-E. Bj\"{o}rk \cite{erik},
which we will describe in this section.

Let \La\ be a (not necessarily commutative) Noetherian ring with
finite injective dimension $d$, i.e $d$ is the minimal integer
 with respect to the property that $ \E^j (M)=0$ for all (left
and right) \La-modules $M$ and integers $j>d$. Of course, this is
equivalent to the condition that both the left and the right
\La-module \La\ has (bounded) injective dimension $d.$ It can be
shown that these left and right injective dimensions are the same
(see \cite{zaks}). The analogous  statement that the left and the
right global homological dimension are the same is a consequence
of the Tor-dimension theorem \cite[4.1.3]{weibel}.\\

 {\em In this section all \La-modules are assumed to be finitely
 generated.}\\

Since  projective \La-modules are reflexive, we get the equality
$$ M=\rhom(\rhom(M,\La),\La)$$ for
 left (or right) \La -modules $M$ in the derived category of
 complexes of \La -modules (more generally, this equality holds for all perfect
complexes). Calculating $ \rhom(\rhom(M,\La),\La)$ by the
bidualizing  complex, the associated filtrations give rise to two
convergent spectral sequences (see \cite{leva} for the
convergence), the first of which degenerates. The second one
becomes

$$ E_2^{p,q}=\E^{p}(\E^{-q}(M)) \Rightarrow \H^{p+q}({\Delta}^{\bullet}(M)
 ), $$
 where $\Delta ^{\bullet}(M)$ is a filtered  complex,
which is exact in all degrees except zero: $\H^0(\Delta
^{\bullet})=M$, i.e.\ there is a canonical filtration $$ \t0 (M)
\subseteq \t1 (M)\subseteq\cdots\subseteq\t{d-1} (M)\subseteq\t{d}
(M)=M$$ on every  module $M$ (The natural numbering of this
filtration arising from the double complex is descending but we
found it more convenient to reverse it).  The convergence of the
spectral sequence implies
$$E^{p,q}_{\infty}=\left\{\begin{array}{cl}
  \t{d-p} (M)/\t{d-p-1}(M) &  \mbox{if } p+q=0, \\
  0 &\mbox{otherwise.}
\end{array}\right.
$$
 (By convention, $\t{i} (M)=0\ \mathrm{for}\ i<0$).

\begin{defn}\begin{enumerate}
\item The number $\delta:=min\{i\mid \t{i} (M)=M\}$ is called the
{\em dimension} $\delta(M)$ of a \La-module $M\neq \{0\}$. We set
$\delta(\{0\})=- \infty.$
\item If $M$ is a  \La-module we say that it has
{\em pure} $\delta$-dimension if \/ $\t{\delta-1}(M)=0,$ i.e.\ the
filtration degenerates to a single term $M$.
\item A \La-module $M$ is called {\em pseudo-null}, if it is at
least of codimension $2,$ i.e.\ if $\delta(M)\leq d-2.$
\end{enumerate}
\end{defn}

By Grothendieck's local duality theorem, this definition coincides
with the Krull dimension of $supp_{\La}(M)$ if \La\ is a
commutative local Noetherian Gorenstein ring, see \cite[Cor.
3.5.11]{bruns}.\\

 First we want to state some basic facts on the
$\delta$-dimension. The functoriality of the spectral sequence implies

\begin{prop}\label{basic}
\begin{enumerate}
\item If  $0\longrightarrow M'\longrightarrow M\longrightarrow
M''\longrightarrow 0$ is an exact sequence of \La-modules then
$$\t{i}(M')\subseteq\t{i}(M) \fa i$$ and
$\delta(M'')\leq\delta(M).$
\item $\t{i}(\bigoplus_k {M_k})=\bigoplus_k\t{i}(M_k)$ and $\delta(\bigoplus_k M_k)=\max_k \delta(M_k).$
\end{enumerate}
\end{prop}

In order to analyze this spectral sequence more closely, the
Auslander condition (for not necessarily commutative rings) is
essential:

\begin{defn}\label{defAuslander}
\begin{enumerate}
\item If $M\neq 0$ is a  \La-module, then
$$j(M):=min\{i\mid\E^i(M)\neq 0\}$$ is called the {\em grade} of
$M.$ By convention,  $j(\{0\})= \infty.$
\item A Noetherian ring \La\ is called Auslander-Gorenstein ring if it has
finite injective dimension and the following Auslander condition
holds: For any \La-module $M$, any integer $m$ and any submodule
$N$ of $\E^m(M),$ the grade of $N$ satisfies $j(N)\geq m.$
\item  A Noetherian ring \La\ is called Auslander regular ring if it has
finite global homological  dimension and the  Auslander condition
holds.
\end{enumerate}
\end{defn}

\begin{rem} Let \La\ be a commutative ring. Then \La\ is Auslander-\linebreak Gorenstein if and only
if it is Gorenstein (in the usual sense). Similarly, \La\ is
Auslander regular if and only if it is regular (in the usual
sense) and of finite Krull dimension. (The implications concerning
the injective, respectively global homological dimensions are well
known, for the Auslander condition see \cite[Cor. 4.6,Prop.\
4.21]{aus} )
\end{rem}
 In the next section we will prove that $\La=\La(G)$ is Auslander regular for any $p$-adic
Lie group without $p$-torsion. Generally, for this kind of rings
we get the following  properties:

\begin{prop}\label{basic1}
Let \La\ be an Auslander regular ring and $M$ a \La-module. Then
\begin{enumerate}
\item\begin{enumerate}

\item \label{clear} For all $i$, there is an exact sequence of \La-modules\\
 \xymatrix{
  0\ar[r] & {\t{i} (M)/\t{i-1} (M)\ar[r]}  & {\E^{d-i}\E^{d-i}(M)\ar[r]} & Q_i(M)\ar[r] & 0,}

   where $Q_i(M)$ is a subquotient of $\bigoplus_{k\geq1}\E^{d-i+k+1}\E^{d-i+k}(M).$
\item $\t0 (M)=\E^d\E^d(M)$ and $\t1 (M)/\t0 (M)=\E^{d-1}\E^{d-1} (M).$
\item $\t{i} (M)/\t{i-1} (M)=0$ if and only if
$\E^{d-i}\E^{d-i}(M)=0.$

      \end{enumerate}

\item \label{dim-formel}\ $\delta(M)+j(M)=d$ for $M\neq 0.$

\item\begin{enumerate}
\item $j(\E^{i} (M))\geq i$, i.e.\ $\E^j\E^{i} (M)=0\fa j<i.$
\item $\delta (\E^{i}(M))\leq d-i.$
\item  $\E^{j(M)}(M)$ has pure $\delta$-dimension $\delta(M).$
\end{enumerate}

\item \ $\E^{k+j(M)+1}\E^{k+j(M)}\E^{j(M)}(M)=0$ for all $k\geq 1.$
\item \label{puredimen}\begin{enumerate}
\item For all $0\leq i\leq d$, $\E^i\E^i(M)$ is either zero or of pure $\delta$-dimension $d-i.$
\item $M$ has pure $\delta$-dimension if and only if
$\E^i\E^i(M)=0$ for all $i>j(M).$
\end{enumerate}
\item\begin{enumerate}
\item $\delta (\t{i}(M))\leq i.$
\item $\t{i}(M)$ is the maximal submodule of $M$ with
$\delta$-dimension less or equal to $i.$
\item The functor $\t{i} $ is left exact.
\item $\t{i}(M/\t{i} (M))=0.$
\end{enumerate}

\item \ If\ $\delta(M)=0$ then $M$ has finite length.
 \end{enumerate}
\end{prop}

\begin{proof}
Except for (i) (a), (i) (b) and (vi), these properties are all proved in \cite{erik} or
trivial: 
 Prop.\ 1.21, 1.16, Prop.\ 1.18, Remark before 1.19, Cor.\ 1.20,
 Cor.\
1.22 and 1.27., while (i)(a) is proved in \cite[Cor.\ 4.3]{leva}

 The assertion
(i)(b) is clear, as $E_{\infty}^{i,-i}=E_{2}^{i,-i}$ because of
(iii)(a). So let us prove (vi): By (iii), (a) is equivalent to
$j(\t{i} (M))\geq d-i$ and this is  true because of the Auslander
condition using induction (cf.\ the proof of (iii)). The
assertion(b) follows by prop.\ \ref{basic}: If $N$ is a submodule
of dimension $\delta(N)\leq i,$ then $N=\t{i}(N)\subseteq \t{i}
(M).$
\\ Noting prop.\ \ref{basic} (i), we only have to show $N\cap\t{i}
(M)\subseteq\t{i} (N)$ in order to prove left exactness. Since the
first module has dimension $\delta (N\cap\t{i} (M))\leq i,$ this
is a consequence of (b).\\ By (c) the exact sequence
$$0\rightarrow\t{i+1} (M)/\t{i} (M)\rightarrow\t{i+2} (M)/\t{i}
(M)\rightarrow\t{i+2} (M)/\t{i+1} (M)\rightarrow 0$$ induces  the
exact sequence $$0\rightarrow\t{i} (\t{i+1} (M)/\t{i}
(M))\rightarrow\t{i} (\t{i+2} (M)/\t{i} (M))\rightarrow\t{i}
(\t{i+2} (M)/\t{i+1} (M)).$$ The first and third term are zero by
(i) and (iii) as above. Hence (d) follows by induction.
\end{proof}

Assuming the Auslander-condition, prop.\ \ref{basic} can be
sharpened as follows  \cite[Prop.\ 1.8]{erik2}:
\begin{prop}
Let  $0\longrightarrow M'\longrightarrow M\longrightarrow
M''\longrightarrow 0$ be an exact sequence of \La-modules.
\begin{enumerate}
\item If \La\ is Auslander-Gorenstein, then
\[j(M)=\min\{j(M'),j(M'')\}\] holds.
\item If \La\ is Auslander regular, then
\[\delta(M)=\max\{\delta(M'),\delta(M'')\}\] holds.
\end{enumerate}
\end{prop}

For the second assertion we used prop.\ \ref{basic1} (ii).

\begin{rem}\label{torchar}
(i)  Using the methods of \cite{fossum}, proposition 6, one can
show the existence of the following exact sequences:
 $$\xymatrix@1@C=12pt{
    {\ 0 \ar[r] } &  {\ \E^{i+1}\du\Omega^i\t{d-i}(M) \ar[r] } &  {\ \t{d-i}(M) \ar[r] } &  {\ \E^i\E^i(M)  \ar[r] } &  {\
    \E^{i+2}\du\Omega^i\t{d-i}(M)\ar[r] }& 0. \
 }$$
Hence, $\t i(M)$ can also be obtained recursively by the
  formula\linebreak
$\t{d-i-1}(M)=\E^{i+1}\du\Omega^i\t{d-i}(M)$ and similarly, we get
a description for \linebreak
$Q_{d-i}(M)\cong\E^{i+2}\du\Omega^i\t{d-i}(M).$ The same arguments
yield for a $\Lambda$-module $M$ with $j(M)\geq j$ the
isomorphisms $$\E^{j+k}\E^j(M)\cong\E^{j+k+2}\du\Omega^j(M)\mbox{
for }k\geq 1.$$ (ii) In particular, $\t{d-1} (M)=\E^1\du
(M)=\tor_\La M,$ i.e.\ the torsion submodule of $M$ is the maximal
submodule of codimension greater or equal than  $1.$  That means
that $M$ is \La-torsion if and only if it is at least of
codimension $1,$ and \La-torsion-free if and only if $M$ is of
pure dimension $d.$
\end{rem}

 The following
class of \La-modules satisfies some duality relations:

\begin{defn} A \La-module $M\neq 0$ is called Cohen-Macaulay
 if \linebreak $j(M)=\pd_{\La}(M)$ holds, i.e.\ if $\E^i (M) =0 \fa i\neq j(M)$.
\end{defn}

\begin{prop}Let $\Lambda$ be an Auslander regular ring.
\begin{enumerate}
\item Let $M$ be a Cohen-Macaulay module of dimension $j$. Then $$\E^j\E^j (M)=M.$$
\vspace{-0.5cm}
\item In particular, if $\delta (M)=0$, then $$\E^d\E^d (M)=M.$$
\end{enumerate}
\end{prop}

\begin{proof}
In both cases the spectral sequence degenerates.
\end{proof}

One could hope that any $\Lambda$-module $M$ can be decomposed
into Cohen-Macaulay modules in the following sense: there is an
filtration of $M$ such that  the $i$th subquotient is
Cohen-Macaulay of dimension $i.$ But it is easy to construct
counterexamples which show that in general such a filtration does
not exist. Nevertheless, there is a different type of filtration:
 Auslander and Bridger
proved the existence of a {\em spherical filtration} (up to
homotopy, i.e.\ after adding a projective summand $P$)
 $$M_d\subseteq M_{d-1}\subseteq\cdots\subseteq M_1\subseteq M_0=M\oplus
 P,$$
 the subquotients of which form {\em spherical} or {\em
Eilenberg-MacLane} modules of type $\E^i(M),$ i.e.\ for $1\leq
i\leq d$
 $$\E^j(M_{i-1}/M_i)\cong\left\{\begin{array}{ll}
   \E^i(M) & \mbox{if }j=i  \\
   0 & \mbox{if }j\neq i,0 \
 \end{array}\right..$$
Fossum \cite{fossum} compared the spherical filtration to the
filtration $\t i(M)$ for a commutative Gorenstein ring and proved
(\cite{fossum}, prop.\ 9) that their ``torsion parts" agree for
$i<d$
 \begin{eqnarray*}
\t i(M)&\cong&\t{d-1}(M_{d-i-1})\\
       &\cong&\t{i}(M_{k})\mbox{ for all } k<d-i.
 \end{eqnarray*}
The proof generalizes at once to the non-commutative case.

\begin{prop}\label{pd}
Let $\Lambda$ be an Auslander regular ring.
 A \La-module $M$ with projective dimension $\pd_{\La}(M)=k$ has no non-trivial
  submodule of dimension less or equal to $d-k-1$, i.e.\ $\t{d-k-1} (M)=0.$
\end{prop}
\begin{proof}
See prop.\ \ref{basic1}, (i) (b).
\end{proof}

The next result extends a well known result for  commutative
regular rings (see for example \cite{nsw}, cor. 5.1.3).

\begin{prop}\label{reflexiv}
Let $\Lambda$ be an Auslander regular ring. \begin{enumerate}
\item For any $\Lambda$-module $M$, the module $\E^0(M)$ is
reflexive: $$\E^0(M)\cong\E^0\E^0\E^0(M).$$\vspace{-0,5cm}
\item Assume that $d\geq 2$ and $\delta(M)=d-i.$ Then
 $$\E^i(M)\cong\E^i\E^i\E^i(M).$$
\end{enumerate}

\end{prop}

\begin{proof}
Let $N:=\E^0(M)$ and  apply proposition \ref{basic1} (iv) to
conclude that \linebreak
$\bigoplus_{k\geq1}\E^{k+1}\E^{k}(N)=\bigoplus_{k\geq1}\E^{k+1}\E^{k}\E^0(M)=0,$
i.e.\ $Q_d(N)=0.$ Since we already know by (iii)(c) that $N$ is of
pure dimension $d$ (if $N\neq 0$), the statement (i) follows
considering (i)(a). The proof of (ii) is  analogous.
\end{proof}

\begin{cor}
For any $i$ it holds
\begin{enumerate}
\item $\E^i\E^i\E^i\E^i(M)\cong\E^i\E^i(M)$ and
\item $\E^i\E^i\t{d-i}(M)\cong \E^i\E^i(M).$
\end{enumerate}

\end{cor}

\begin{proof}
To prove (i) assume first that $\delta(\E^i(M)=d-i$. Applying the
proposition to the module $\E^i(M)$ gives the result while  in the
second case, i.e.\ $j(\E^i(M)>i,$ the module $\E^i\E^i(M)$ is zero
anyway. Noting that $j(Q_i(M))\geq i+2,$ the second assertion
follows at once calculating the long exact $\E^i$-sequence of
 $$\xymatrix@1{
    {\ 0 \ar[r] } &  {\ \t{d-i-1}(M)/\t{d-i}(M) \ar[r] } &  {\ \E^i\E^i(M) \ar[r] } &  {\ Q_i \ar[r] } &  {\
    0. } \
 }$$
\end{proof}

\subsection{Modules up to pseudo-isomorphisms}
As in the commutative case we say that a homomorphism
$\varphi:M\to N$ of \La-modules is a {\em pseudo-isomorphism} if
its kernel and cokernel are pseudo-null. A module $M$ is by
definition pseudo-isomorphic to a module $N,$ denoted
 $$M\sim N,$$
if and only if there exists a pseudo-isomorphism from $M$ to $N.$
In general, $\sim$ is not symmetric even in the $\zp$-case (cf.\
\cite[V\S3,ex.1]{nsw}). While in the commutative case $\sim$ is
symmetric at least for torsion modules (see the first remark after
prop.\ 5.17 in \cite{nsw}), we do not know whether this property
still holds in the general case. \\ If we want to reverse
pseudo-isomorphisms, we have to consider the quotient category
$\mod/{\mathcal{PN}}$ with respect to subcategory $\mathcal{PN}$
of pseudo-null \La-modules, which is a Serre subcategory, i.e.\
closed under subobjects, quotients and extensions. By definition,
this quotient category is the localization $({\mathcal{
PI}})^{-1}\mod$ of \mod\ with respect to the multiplicative system
${\mathcal{ PI}}$ consisting of all pseudo-isomorphisms (see
\cite[ex.\ 10.3.2]{weibel}). Since \mod\ is well-powered, i.e.\
the family of submodules of any module $M\in\mod$ forms a
set,
these localization exists, is an abelian category and the
universal functor $q:\mod\to\mod/{\mathcal{PN}}$ is exact (see
\cite[p.\ 44ff]{swan}). Furthermore, $q(M)=0$ in
$\mod/{\mathcal{PN}}$ if and only if $M\in\mathcal{PN}.$

Note that for any pseudo-isomorphism $f:M\to N$ the induced
homomorphism $\E^1(f)$ is an pseudo-isomorphism (with respect to
an analogue subcategory, also denoted by $\mathcal{PN},$ of
$\mathrm{mod-}\La$), too. By the universal property of the
localization, we obtain a functor
\[ E^1(-):\mod/{\mathcal{PN}}\to \mathrm{mod-}\La/{\mathcal{PN}},\]
which is exact if it is restricted to the full subcategory
$\La\mathrm{-mod}^{\geq 1}/{\mathcal{PN}}$ of \linebreak
$\mod/{\mathcal{PN}}$ consisting of all \La-modules of codimension
greater or equal to $1.$ The next result shows that
\[\E^1\circ\E^1\cong id:\La\mathrm{-mod}^{\geq 1}/{\mathcal{PN}}\to \La\mathrm{-mod}^{\geq 1}/{\mathcal{PN}}\]
is a natural isomorphism of functors:

\begin{prop}\label{E^1-pseudo} Let \La\ be an Auslander regular ring.
\begin{enumerate}
\item Any torsion-free module $M$ embeds into a reflexive module
with pseudo-null cokernel.
\item Any torsion module $M$ is pseudo-isomorphic to
$\E^1\E^1(M).$
\end{enumerate}

\end{prop}

\begin{proof}
Observe that $\t{d-1}(M)$ is the maximal torsion submodule in this
case. Hence, the exact sequence in (i) (a) for $i=d$ respectively
$i=d-1$ proves both statements taking under consideration (iii)(b)
and proposition \ref{reflexiv}.
\end{proof}

\begin{cor}
A homomorphism $f:M\to N$ of \La-torsion modules is a
pseudo-isomorphism if and only if the induced homomorphism
$\E^1(f):\E^1(N)\to\E^1(M)$ is.
\end{cor}

Since \La\ is Noetherian it follows readily that any object in
$\La\mathrm{-mod}^{\geq 1}/{\mathcal{PN}}$ has the ascending chain
condition (A.C.C.) (see \cite[II.ch.2]{swan}). But using the
natural isomorphism $\E^1\circ\E^1\cong id$ its immediate that
also the descending chain condition (D.C.C.) holds in this
category. A consequence of this observation is that any object in
$\La-\mathrm{mod}^{\geq 1}/{\mathcal{PN}}$ has finite length.
Moreover, the Krull-Schmidt-Theorem holds (loc. cit. thm. 2.2):

\begin{prop}\label{krull-schmidt} Let \La\ be an Auslander regular ring. Then
any  object \linebreak $q(M)\in\La\mathrm{-mod}^{\geq
1}/{\mathcal{PN}}$ decomposes uniquely (up to isomorphism) into a
finite direct sum of indecomposable objects $q(M_i):$
\[M\equiv\bigoplus_i M_i\ \mathrm{ mod }\ \mathcal{PN}.\]
\end{prop}

\begin{prop}\label{simple} Let \La\ be an Auslander regular ring. Then the
following holds.
\begin{enumerate}
\item The simple objects of $\La\mathrm-{mod}^{\geq
1}/{\mathcal{PN}}$ are cyclic, i.e. of the form $q(\La/I)$ for
some left ideal $I$ of \La.
\item There is an isomorphism of abelian groups \[K_0(\La\mathrm{-mod}^{\geq
1}/{\mathcal{PN}})\cong\bigoplus_{S\in\mathcal{I}}\mathbb{Z}[S]\]
where $\mathcal{I}$ denotes the set of isomorphism classes of
simple objects of \linebreak $\La\mathrm-{mod}^{\geq
1}/{\mathcal{PN}}.$
\end{enumerate}
\end{prop}

The class $[M]$ of a $1$-codimensional \La-module  $M$ in
$K_0(\La\mathrm{-mod}^{\geq 1}/{\mathcal{PN}})$ should be thought
of as "generalized characteristic ideal." At least if \La\ is a
commutative regular ring, this class bears the same information as
the characteristic ideal which is associated with $M$ via the
structure theory of \La-modules up to pseudo-isomorphism (see
\cite[5.1.7,5.18]{nsw}).

\begin{proof}
Let $M$ be a non-pseudo-null \La-module such that $q(M)$ is a
simple object. Then there exists a  $m\in M$ such that $\La m\cong
\La/ann_\La(m)$ is not pseudo-null, either. Consequently, $0\neq
q(\La m)\subseteq q(M)$ and by the simplicity of $q(M)$ this
inclusion cannot be strict. Taking $I=ann_\La(m),$ this proves (i)
while (ii) is just \cite[II ch.1, thm 1.7]{swan}.
\end{proof}

 Following the structure theory for modules over a {\em
commutative} regular local ring (see \cite[5.1.7,5.18]{nsw}), it
is natural to hope that  the following question has an affirmative
answer

\begin{ques}
Let \La\ be an Auslander regular ring and $M\in\mod.$ Does there
exist an isomorphism in $\mod/{\mathcal{PN}}$
 $$M\cong \tor_\La M\oplus R\mbox{ mod }{\mathcal{PN}},$$
where $R\cong M/\tor_\La M\mbox{ mod }{\mathcal{PN}}$ is a
reflexive \La-module?
\end{ques}

At least for the Iwasawa algebra $\La(G)$ of an extra-powerful
(for the definition see the next subsection) and uniform
pro-$p$-group this should follow from the techniques used by
Coates-Sujatha (\cite{co-su1}) to prove the structure theorem for
torsion modules.

\begin{prop}\label{pseudoisoE1} Let \La\ be an Auslander regular ring. For any \La-module $M$ it holds:
$$\E^1(M)\sim\E^1(\tor_\La M).$$
\end{prop}

\begin{proof}
>From the long exact Ext-sequence we get the exact sequence
 $$\xymatrix@1{
    {\ \E^1(M/\tor_\La M) \ar[r] } &  {\ \E^1(M) \ar[r] } &  {\ \E^1(\tor_\La M) \ar[r] } &  {\ \E^2(X/\tor_\La M).  }    \
    }$$
While the  module on the right hand side is obviously pseudo-null
the first one is so by the following argument: the long exact
Ext-sequence of
 $$\xymatrix@1{
    {\ 0 \ar[r] } &  {\ M/\tor_\La M \ar[r] } &  {\ \E^0\E^0(M) \ar[r] } &  {\ \E^2\du(M) \ar[r] } &  {\
    0 } \
 }$$
tells us that $\E^1(M/\tor_\La M)$ fits into the exact sequence
 $$\xymatrix@1{
    {\ \E^1\E^0\E^0(M) \ar[r] } &  {\ \E^1(M/\tor_\La M) \ar[r] } &  {\ \E^2\E^2\du(M),  }  \
    }$$
i.e.\ it suffices to show that $\E^1\E^0\E^0(M)$ is pseudo-null.
But $\E^1\E^1\E^0\E^0(M)=0$ by \ref{basic1},\ref{puredimen}, (and
$\E^0\E^1\E^0\E^0(M)=0$ anyway), i.e.\ $j(\E^1\E^0\E^0(M))\geq 2$
respectively $\delta(\E^1\E^0\E^0(M))\leq d-2.$
\end{proof}
\POP
\PUSH{graduation.tex}%
\subsection{The graded ring $gr($\La$)$}

An important method to verify the Auslander condition of a ring
\La\ consists of endowing \La\ with a suitable filtration and
studying the associated graded ring $gr(\La).$ By a filtration  on
a ring \La\ we mean an increasing (!) sequence of additive
subgroups $\Sigma_{i-1}\subseteq\Sigma_{i}\subseteq\Sigma_{i+1}$
satisfying $\bigcup\Sigma_{i}=\La$ and $\bigcap\Sigma_{i}=0$ and
the inclusions $\Sigma_{i}\Sigma_{k}\subseteq\Sigma_{i+k}$ hold
for all pairs of integers $i$ and $k$. The main example  on a
local ring is the \M-adic filtration with $\Sigma_{-i}=\M^i \fa
i\geq 0$ (by convention, $\M^0 =\La$ ). For our aim the {\em
closure condition} will be crucial:

\begin{defn}
The filtration $\Sigma$ satisfies the {\em closure condition} if
the additive subgroups $\Sigma_{m_1} u_1+\cdots +\Sigma_{m_{s}}
u_s$ and $u_1\Sigma_{m_1}+\cdots +u_s \Sigma_{m_{s}}$ are closed
with respect to the topology induced by $\Sigma$ for any finite
subset $u_1 ,\ldots ,u_s$ in \La\ and all integers $m_1 ,\ldots
,m_s .$
\end{defn}

\begin{lem}
Let $G$ be a $p$-adic analytic pro-$p$-group. Then the \M-adic filtration on $\La (G)$
satisfies the closure condition.
\end{lem}

\begin{proof}
Note that the $\frak M$-adic topology on \La\ coincides with the
$(\m ,I)$-topology (cf.\ \cite[(5.2.15)]{nsw}). Since \M\ is a
two-sided ideal of \La\ the subgroup $\M^{i-m_1} u_1+\cdots
+\M^{i-m_{s}} u_s$, $u_1\M^{i-m_1}+\cdots +u_s \M^{i-m_{s}}$ is a
finitely generated left, right ideal, respectively. Hence, these
subgroups are compact as continuous images of the compact module
$\La ^n$ for some $n$.
\end{proof}

Put $gr(\La)=\bigoplus\Sigma_i /\Sigma_{i-1}$, which is called the
associated graded ring of \La\ with respect to the filtration
$\Sigma .$ The above lemma admits applying the following theorem
of Bj\"{o}rk to certain completed group rings:

\begin{thm}[Bj\"{o}rk]
\begin{enumerate}
\item Assume that $gr(\La)$ is a  Auslander regular ring and that $\Sigma$ satisfies the
closure condition. Then \La\ is a Auslander regular ring.
\item In the situation of {\rm (i),} the equality $j(M)=j(gr(M))$ holds.
 If, in addition, $gr(\La)$ is commutative and of pure dimension $d$, then also $\delta(M)=\dim (gr(M))$ holds,
where $\dim (gr(M))=\dim (supp_{gr(\La)}(gr(M)))$ is the Krull
dimension of $gr(M)$.
\end{enumerate}
\end{thm}

\begin{proof}
See \cite[Theorems 4.1,4.3]{erik} and also \cite[Thm. 3.9. and
Remark]{erik2}. For the last equality note that
\begin{eqnarray*}
\dim (gr(M))&=&\max \{\dim(gr(M)_{\p}\mid \p\ maximal\ ideal\ of\
gr(\La)\}\\ &=&d-\min\{j(gr(M)_{\p})\mid \p\ maximal\ ideal\ of\
gr(\La)\}\ \\ &=&d- j(gr(M)\\ &=&d-j(M)\\ &=&\delta(M),
\end{eqnarray*}
where we used prop.\ \ref{basic1}, (ii), and the fact that
localization commutes with Ext-groups, if all objects are
Noetherian.
\end{proof}

Our task will be to determine the structure of $gr(\La (G))$.
Before stating the next, result we recall that a pro-$p$-group $G$
is called extra-powerful, if the relation $[G,G]\subseteq G^{p^2}$
holds. Furthermore, note that $gr(\zp)\cong\fp[X_0 ]$ if \zp\ is
endowed with the \m-adic filtration.

\begin{thm}\label{gr}
Let $G$ be a uniform and extra-powerful pro-$p$-group of dimension
$\dim(G)=r.$ Then there is a $gr(\zp )$-algebra-isomorphism
$$gr(\La (G))\cong\fp[X_0,\ldots ,X_r ].$$ In particular, $gr(\La
(G))$ is a commutative regular Noetherian ring.
\end{thm}

A consequence of  Lazard's results is the

\begin{rem}
Any compact $p$-adic analytic group contains an open
characteristic subgroup, which is an uniform and extrapowerful
pro-$p$-group (cf.\ \cite[Cor. 9.36]{dsms} and \cite[Prop.\
8.5.3]{wi})
\end{rem}

For the proof of the theorem we need some more terminology. Let
$G$ be an uniform pro-$p$-group with a minimal system of
(topological) generators $\{x_1 ,\ldots ,x_r \}$, in particular
$\dim(G)=r.$ Then the lower $p$-series is given by $P_1(G)=G,\
P_{i+1}(G)=(P_i (G))^p,\ i\geq 1.$ This filtration defines a
$p$-valuation $\omega :G\longrightarrow
\n_{>0}\cup\{\infty\}\subseteq\mathbb{R}_{>0}\cup\{\infty\}$ of
$G$ in the sense of Lazard via $\omega(g):=\sup\{i\mid g\in P_i
(G) \}$, which induces a filtration on the group algebra $\zp [G]$
of the underlying abstract group of $G,$ too (cf.\ \cite[Chap.\
III, 2.3.1.2]{la}).

\begin{lem}
The filtration on $\mathbb{Z}_p[G]$, induced by $\omega$, is the
$\M_d$-adic one, where $\M_d =\m + I_d (G)$ with the augmentation
ideal $I_d(G)$ of $\zp [G].$
\end{lem}

\begin{proof}
Conferring  the proof of Lemma III, (2.3.6) in \cite{la}, the
induced filtration is given by the following ideals in $\zp [G]$,
$n\in \n:$ $A_n$ is generated as $\zp$-module by the elements
$p^l(g_1-1)\cdots (g_m -1)$ where $l,m\in \n,\ g_i\in G$ and
$l+\omega (g_1)+\ldots +\omega(g_m)\geq n$, whereas the
$\M_d$-adic filtration is defined by the ideals $\M_d^n$, which
are generated (over $\zp [G]$) by the elements $p^l(g_1-1)\cdots
(g_m-1)$, where $l,m\in\n,\ g_i\in G$ and $l+m=n$. Since
$\omega(g)\geq1\fa g\in G$ the ideal $\M_d^n$ is contained in
$A_n$.  The converse is a consequence of the following \\[-0.2cm]

 {\em Claim:} Let $g\in G$ with $\omega(g)=t\geq
1,$ then $g-1\in\M_d^t.$\\

\vspace*{-0.3cm} Since $G$ is uniform, the map $G\longrightarrow
P_t(G),$ which assigns $g^{p^{t-1}}$ to $g,$ is surjective (cf.\
\cite[lemma 4.10]{dsms}), i.e.\ there exists an element $h\in G$
with $g=h^{p^{t-1}}$. Writing
$$g-1=(1+(h-1))^{p^{t-1}}-1=\sum_{k\geq 1} \left(
\begin{array}{c}
  p^{t-1} \\
  k
\end{array}\right)(h-1)^k,$$ one verifies that $g-1\in \M^t_d$, because
$v_p(\left(\begin{array}{c}
  p^{t-1} \\
  k
\end{array}\right))=t-1-v_p(k)\geq t-k,$  i.e.\ $\left(\begin{array}{c}
  p^{t-1} \\
  k
\end{array}\right)(h-1)^k\in \M^t_d$.
\end{proof}

\begin{lem}\label{madic}
The $\M_d$-adic filtration on $\zp [G]$ induces  the $\M$-adic
filtration on $\zp\kl G\kr $.
\end{lem}

\begin{proof}
The ideals defining the induced filtration are just the closure
$\overline{\M_d^n}$ of $\M_d^n\subseteq\zp[G]\subseteq\zp\kl G\kr
$ with respect to the $\M$-adic topology on $\zp\kl G\kr $. Since
they
 contain the elements $p^l(g_1-1)\cdots (g_m-1)$ with $\
l,m\in\n,$ $g_i\in\{x_1,\ldots ,x_r\}$ and $l+m=n$, which generate
$\M^n$ as ideal of $\zp\kl G\kr $, they contain  $\M^n$, too. On
the other hand $\M^n$ is closed and contains all the generators of
the $\zp[G]$-ideal $\M^n_d:$ $p^l(g_1-1)\ldots (g_m-1), \
l,m\in\n,\ g_i\in G$. This proves the lemma.
\end{proof}

\vspace*{-0.3cm} Now we can prove theorem \ref{gr}.

\begin{proof}
Since $gr(G)=\bigoplus P_i(G)/P_{i+1}(G)$ is a Lie algebra, which
is free of rank $r$ as $gr(\zp)$-module, we get the following
inclusion: \begin{eqnarray*} gr(G)\subseteq
Ugr(G)&\cong&gr(\zp[G])\\ &\cong&gr(\zp\kl G\kr ),
\end{eqnarray*} where the first equation holds by \cite[Chap.\ III,
2.3.3]{la} and $Ugr(G)$ is the enveloping algebra of the Lie
algebra $gr(G)$, whereas the second one is a consequence of lemma
\ref{madic}. According to \cite[Theorem 8.7.7]{wi}, the graded
ring $gr(\zp\kl G\kr )$ is commutative ($G$ is assumed to be
extra-powerful), i.e.\ $$Ugr(G)\cong
gr(\zp)[X_1,\ldots,X_r]\cong\fp[X_0,\ldots,X_r].$$
\end{proof}

\vspace*{-0.35cm} As an  important consequence we obtain the

\begin{thm}\label{Lambda-Auslanderreg}
Let $G$ be a compact $p$-adic analytic group without $p$-torsion. Then the completed
group ring $\La(G)$ is an Auslander regular ring.
\end{thm}

\begin{proof}
$G$ posses an open characteristic subgroup $N$ which is an
uniform, extra-powerful pro-$p$-group. By the theorem of Bj\"{o}rk and
theorem \ref{gr}, $\La(N)$ is an Auslander regular ring, because
$gr(\zp\kl N\kr )$ has this property as a regular commutative
Noetherian ring (cf.\ \cite[pp.\ 65-69]{erik1}). But
$\E^i_{\La(G)}(M)\cong\E^i_{\La(N)}(M)$ as $\La(N)$-modules for
any $\La(G)$-module $M$, by which the Auslander condition is
easily verified.
\end{proof}

If one is not  interested in the $\zp$-torsion submodule of a
\La-module $M$ it might be convenient to replace $M$ by its
maximal $\zp$-torsion-free quotient $M/\tor_{\mathbb{Z}_p} M.$
This leads to the study of the subcategory
$\La\mathrm{-mod}^{\mathrm{fl}}$ of $\La\mathrm{-mod}$ consisting
of $\zp$-torsionfree or equivalently $\zp$-flat \La-modules. This
category is closely related to the category $\qp\kl G
\kr\mathrm{-mod}$ of finitely generated $\qp\kl G \kr$-modules
where \[{{\Bbb Q}_p}\kl G \kr:=\zp\kl G \kr \otimes_{\mathrm{Z}_p} \qp\] by
definition. Note that the latter ring is  Noetherian whenever
$\zp\kl G\kr$ is. We let $\La\mathrm{-mod}^{\mathrm{fl}}_\qp$
denote the category with the same objects as
$\La\mathrm{-mod}^{\mathrm{fl}}$ and such that
 \[\Hom_{\La\mathrm{-mod}^{\mathrm{fl}}_\qp}(M,N):=\Hom_\La(M,N)\otimes_\mathbb{Z}\mathbb{Q}\]
for any two \La-modules $M,N$ in $\La\mathrm{-mod}^{\mathrm{fl}}$
with the composition of homomorphisms being the
$\mathbb{Q}$-linear extension of the composition in
$\La\mathrm{-mod}^{\mathrm{fl}}.$ It is clear that the functor $M
\mapsto M_\qp:=M\otimes_\zp\qp$ induces an equivalence of
categories
 \[\Lambda\mathrm{-mod}^{\mathrm{fl}}_\qp\cong\qp\kl G
\kr\mathrm{-mod}.\]

\begin{rem}
Schneider and Teitelbaum proved in \cite{sch-teit} that, for a
compact $p$-adic Lie group, these categories are anti-equivalent
to the category of certain Banach space representations of
$G.$\end{rem}

The following proposition, which is standard, allows to calculate
Ext-groups in $\qp\kl G \kr\mathrm{-mod}$ via Ext-groups in
$\mod.$

\begin{prop}
For any $M,N$ in \mod  we have canonical isomorphisms
 \[\mathrm{Ext}_\La^i(M,N)\otimes_\zp\qp\cong\mathrm{Ext}_{\qp\kl G \kr}^i(M_\qp,N_\qp) \fa i\geq0.\]
\end{prop}

Now, we are able to derive the Auslander regularity of  $\qp\kl G
\kr$ from the Auslander regularity of $\La(G).$

\begin{thm}\label{qp-Auslander}
Let $G$ be a compact $p$-adic analytic group without $p$-torsion.
Then the  group ring $\qp\kl G \kr$ is an Auslander regular ring
of (projective) dimension $\pd \qp\kl G \kr =\cd_p(G).$
\end{thm}

\begin{proof}
We set $r:=\cd_p G.$ Let $M$ be any finitely generated $\qp\kl G
\kr$-module and choose $N\in\La\mathrm{-mod}^{\mathrm{fl}}$ with
$M\cong N_\qp.$ Since $N$ is \zp-torsionfree, we have
$E^{r+1}_\La(N)=0$ by theorem \ref{ja-tor} (iii), hence $\pd_\La
N\leq d$ using (\ref{pdviaExt}). By the previous proposition,
$\E^{r+1}_{\qp\kl G \kr}(M)$ vanishes, too, i.e. $\pd_{\qp\kl G
\kr} M\leq r.$ On the other hand the projective dimension of $\qp$
is $r,$ because $\E^r_{\qp\kl G
\kr}(\qp)=\E^r(\zp)\otimes_\zp\qp\cong\qp$ by \cite[2.6]{ja-is}.
It follows that $\pd \qp\kl G \kr =r.$

Now we will verify the Auslander condition: Since any $\qp\kl G
\kr$-submodule of $\E^i_{\qp\kl G
\kr}(M)\cong\E^i_\La(N)\otimes_\zp\qp$ has the form $L_\qp$ for
some \La-submodule $L\subseteq\E^i_\La(N)$ we see that
 \[\E^j_{\qp\kl G
\kr}(L_\qp)\cong\E^j_\La(L)\otimes_\zp\qp=0,\ j<i,\] because \La\
is Auslander regular by theorem \ref{Lambda-Auslanderreg}.
\end{proof}

\subsection{The $\mu$-invariant}

For the purpose to study the $p$-torsion part $\tor_{\mathbb{Z}_p}
M$ of a \La-module $M$ we are also interested in the rings
$\z/p^m\kl G\kr \cong\La(G)/p^m,$ especially the ring $\fp\kl G\kr
,$ and will consider the change of rings
$\La(G)\rightarrow\La/p^m.$ For a $\La/p^m$-module $M$ there
exists a convergent spectral sequence (see \cite[Ex.
5.6.3]{weibel})
 $$\Ext^i_{\La/p^m}(M,\Ext^j_\La(\La/p^m,\La))\Rightarrow\Ext_\La^{i+j}(M,\La).$$
We should mention that here $\Ext^j_\La(\La/p^m,\La)$ is a left
$\La$- and $\La/p^m$-module by functoriality and the right
\La-structure of the bi-module $\La/p^m.$ Using the free
resolution
 $$\xymatrix@1{
    {\ 0 \ar[r] } &  {\ \La \ar[r]^{p^m} } &  {\ \La \ar[r] } &  {\ \La/p^m \ar[r] } &  {\ 0, } \
 }$$
it is easy to calculate that

 $$\Ext^j_{\La}(\La/p^m,\La)\cong\left\{\begin{array}{cl}
   \Lambda/p^m & \mbox{if } j=1, \\
   0 & \mbox{otherwise.} \
 \end{array}\right.$$
Hence the spectral sequence degenerates to
 $$\E^i_{\La/p^m}(M)\cong\E^{i+1}_\La(M)$$
for any $\La/p^m$-module $M$ and any integer $i.$ We obtain the
following

\begin{thm}\label{fp-Auslander}
Let $G$ be a compact $p$-adic analytic group without $p$-torsion
and $m$ any natural number. Then \begin{enumerate}
\item $\zp/p^m\kl G\kr $ is an Auslander-Gorenstein ring of injective
dimension $\cd_p(G).$
\item $\fp\kl G\kr $ is an Auslander regular ring of  dimension
$\cd_p(G).$
\end{enumerate}
\end{thm}

\begin{rem}
The same arguments prove that $\La(G)/(f)$ is a Auslander
Gorenstein ring of injective dimension $\cd_p(G)$ for any element
$f$ of the center of $\La(G)$ or even more general for any
$f\in\La(G)$ such that the left ideal $(f):=\La(G)f$ is two-sided.
 \end{rem}

\begin{proof}
>From the above formula we derive that $\La/p^m$ has finite
injective dimension $\cd_p(G).$ On the other hand it is well known
that the projective dimension of $\fp\kl G\kr $ is equal to
$\cd_p(G)$ (see \cite[V\S2Ex.5]{nsw}). Hence it suffices to verify
the Auslander condition: For a $\La/p^m$- module $M$ let
$N\subseteq\E^i_{\La/p^m}(M)$ be a $\La/p^m$-submodule which we
will also consider as \La-submodule of $\E^{i+1}_\La(M).$ Applying
again the above isomorphism, we see that
$\E^j_{\La/P^m}(N)\cong\E^{j+1}_\La(N)=0$   for any integer $j<i$
because \La\ fulfills the Auslander condition.
\end{proof}

A different possibility to prove (ii) of the previous theorem
would be to imitate the proof of theorem \ref{Lambda-Auslanderreg}
using the analogue of theorem \ref{gr}: if $G$ is a uniform
pro-$p$-group of dimension $d,$ then there is an isomorphism
 $$gr(\fp \kl G\kr )\cong\fp[X_1,\ldots ,X_r ],$$
where $\fp \kl G\kr $ is endowed with its $\M$-adic filtration
(see \cite[8.7.10]{wi}). In particular, $\fp \kl G\kr $ has no
zero divisors for uniform $G$ (\cite[8.7.9]{wi}).

In order to measure the size of the $p$-torsion part of a
\La-module we have (as usual) the $\mu$-invariant which is defined
as follows.

\begin{defn}
\label{mu} Assume that $G$ is a $p$-adic Lie group without
$p$-torsion such that $\fp\kl G\kr $ is integral. For any
$\La(G)$-module $M$ we define its $\mu$-invariant $\mu(M)$ as
$$\mu(M)=\rk_{\fp\kl G\kr }\bigoplus_{i\geq 0}
{_{p^{i+1}}M}/{_{p^{i}}M},$$ where ${_{p^{0}}M}=0$ by convention.
Observe that the sum is finite because \La\ is Noetherian.
\end{defn}

Susan Howson has defined the $\mu$-invariant in a similar,
equivalent way and she has independently studied its properties in
\cite{howson2000}. In particular, she expresses the
$\mu$-invariant of a \La-module $M$ to the Euler characteristic of
$\tor_\zp M.$

Note that the $\mu$-invariant only depends on the \La- resp.\
$\zp$-torsion-\linebreak submodule: $\mu(M)=\mu(\tor_\La
M)=\mu(\tor_\zp M)=\mu({_{p^m}M})$ for $m$ sufficiently large.
With respect to the vanishing we have the following
characterization:

\begin{rem}\label{mu0}
Since $\xymatrix@1{{\ {_{p^{i+1}}M}/{_{p^{i}}M}
\ar@{^(->}[r]^(.7){p^i} } & {\ _pM }}$ the following is equivalent
 \begin{eqnarray*}
  \mu(M)=0 & \Leftrightarrow &\mu(_pM)=0 \\
  &\Leftrightarrow & {_pM} \mbox{ is $\fp\kl G\kr $-torsion}\\
 &\Leftrightarrow & {_pM} \mbox{ is a pseudo-null $\La$-module.}
\end{eqnarray*}
For the latter equivalence we used again the above isomorphism.
\end{rem}

The next proposition shows that the $\mu$-invariant is in fact an
invariant ``up to pseudo-isomorphism", i.e.\ it factors through
the quotient category $\mod/{\mathcal{PN}}.$ In particular, our
definition of $\mu$ generalizes the usual definition via the
structure theory if $G$ is isomorphic to $\mathbb{Z}_p^r$ for some
$r.$

\begin{prop}
Let $G$ be a $p$-adic analytic group without $p$-torsion such that
both $\La=\La(G)$ and $\La/p$ are integral. Then
 $$M\sim N \mbox{ implies }\mu(M)=\mu(N).$$
\end{prop}

\begin{proof}
The statement will follow if it holds in the two special cases of
exact sequences
\begin{enumerate}
\item[(a)] $\xymatrix@1{
   {\ 0 \ar[r] } &  {\ Q \ar[r] } &  {\ M \ar[r] } &  {\ N \ar[r] } &  {\
   0, }
}$
\item[(b)]$ \xymatrix@1{
   {\ 0 \ar[r] } &  {\ M \ar[r] } &  {\ N \ar[r] } &  {\ Q \ar[r] } &  {\
   0, }
}$
\end{enumerate}
where $Q$ is pseudo-null. More generally, we consider a short
exact sequence of \La-modules
 $$\xymatrix@1{
    {\ 0 \ar[r] } &  {\ X \ar[r] } &  {\ Y \ar[r] } &  {\ Z \ar[r] } &  {\
    0. } \
 }$$
The snake lemma implies the exactness and commutativity of the
following diagram
 $$\xymatrix{
   0\ar[r] & {_{p^n}X}\ar[r]\ar@{^{(}->}[d] & {_{p^n}Y}\ar[r]\ar@{^{(}->}[d] & {_{p^n}Z}\ar[r]\ar@{^{(}->}[d] & X/p^n\ar[d]^p \\
   0\ar[r]  & {_{p^{n+1}}X}\ar[r] & {_{p^{n+1}}Y}\ar[r] &{_{p^{n+1}}Y}\ar[r] & X/p^{n+1} .\
 }$$
Again by the snake lemma we obtain the exact sequences
\begin{eqnarray*}
\xymatrix@1{
   {\ 0 \ar[r] } & {\ {_{p^{n+1}}X}/ {_{p^n}X}\ar[r] } &  {\ {_{p^{n+1}}Y}/ {_{p^n}Y} \ar[r] } &  {\ A_{n+1}/A_n \ar[r] } &    {\
    0, } \
 }&&\\
\xymatrix@1{
  {\ 0 \ar[r] }& {\ K_n \ar[r] } &  {\ A_{n+1}/A_n  \ar[r] } &  {\ {_{p^{n+1}}Z}/ {_{p^n}Z} \ar[r] } &  {\ B_{n+1}/B_n \ar[r] } &  {\
   0, }
}
\end{eqnarray*}
where $A_i$ denotes the image of $_{p^i}Y$ in  $_{p^i}Z$ with
cokernel $B_i,$ the latter module considered as submodule of
$X/p^i,$ and $K_n:=\ker(B_n\to B_{n+1}).$ In case (b)
$A_{n+1}/A_n$ is a pseudo-null \La-module because
$A_{n+1}\subseteq Z.$ Hence $\rk_{\fp\kl G\kr }A_{n+1}/A_n=0$ by
remark \ref{mu0}. In case (a) $\rk_{\fp\kl G\kr }{_{p^{n+1}}X}/
{_{p^n}X}=0$ by the same argument. Furthermore, $K_n\subseteq
X/p^n,$ $B_{n+1}\subseteq X/p^{n+1}$ and finally $B_{n+1}/B_n$ are
pseudo-null, too.
\end{proof}

By $\mod(p)$ we shall write the plain subcategory of $\mod$
consisting of $\zp$-torsion modules while by
${\mathcal{PN}}(p)``={\mathcal{PN}}\cap\mod(p)"$ we denote the
Serre subcategory of $\mod(p)$ the objects of which are
pseudo-null \La-modules. In other words $M$ belongs to ${\mathcal{
PN}}(p)$ if and only if it is a $\La/p^n$-module for an
appropriate $n$ such that $\E^0_{\La/p^n}(M)=0.$ Recall that there
is a canonical exact functor
$q:\mod(p)\to\mod(p)/{\mathcal{PN}}(p).$ For the description of
the $p$-torsion part the following result will be crucial.

\begin{prop}\label{p-simple}
Assume that $G$ is  a $p$-adic analytic group without $p$-torsion
such that both $\La=\La(G)$ and $\La/p$ are integral. Then the
following holds:
\begin{enumerate}
\item $q(\La/p)$ is simple object in $\mod(p)/{\mathcal{PN}}(p),$ i.e.\
does not contain any proper subobject.
\item Every object $A$ in $\mod(p)/{\mathcal{PN}}(p)$ has a finite
composition series
 $$0=A_0\subseteq A_1\subseteq \cdots\subseteq A_{i+1}=A$$
of subobjects $A_j$ of $A$ such that $A_{j+1}/A_j\cong q(\La/p)$
for every $i\geq j\geq 0.$ In particular, $q(\La/p)$ is the unique
simple object of $\mod(p)/{\mathcal{PN}}(p).$
\item Any $q(M)$ in $\mod(p)/{\mathcal{PN}}(p)$ has finite
length equal to $\mu(M).$ Thus, $[q(M)]\mapsto
 \mu(M)$ induces  an isomorphism
 $$K_0(\mod(p)/{\mathcal{PN}}(p))\cong\mathbb{Z}.$$
\end{enumerate}
\end{prop}

We need the following lemma which can be proved literally as
\cite[lem 2.25]{howson} because $\fp\kl G\kr $ is both left and
right Noetherian ring without zero divisors and thus it has a skew
field of fractions.

\begin{lem}\label{fp-lemma} With the assumptions of the proposition
let $M$ be a torsion-free $\La/p$-module of rank
$\rk_{\La/p}(M)=m.$ Then there exist free $\La/p$-modules $F, F'$
such that $F\subseteq M,$ $M\subseteq F'$ and both $M/F$ and
$F'/M$ are $\La/p$-torsion, i.e.\ pseudo-null considered as
\La-module. In particular, for any $\La/p$-module of rank $m$
there is an isomorphism $$q(M)\cong q(\La/p)^m.$$
\end{lem}

\begin{proof}
Let $h:q(M)\hookrightarrow q(\La/p)$ be a monomorphism in the
quotient category. By \cite[I 2.9]{swan} there exists a diagram
 $$\xymatrix{
   M &   & {\La/p} \\
     & M'\ar[lu]^f\ar[ru]_g &   \
 }$$
in $\mod(p)$ with $f$ a pseudo-isomorphism in $\mod$ such that
$$\xymatrix{
   q(M) \ar[rr]^h&   & {q(\La/p)} \\
     & q(M')\ar[lu]^{q(f)}\ar[ru]_{q(g)} &   \
 }$$
commutes. Since $h$ is a monomorphism and $q(f)$ an isomorphism,
$\ker(g)$ must be in ${\mathcal{PN}}(p).$ Since
$M'/\ker(g)\subseteq\La/p,$ we can consider its $\La/p$-rank which
can be either $1$ or $0$. In the first case we conclude that $g$
is a pseudo-isomorphism, i.e.\ $q(g)$ is an isomorphism, while in
the second case $M'/\ker(g)$ and hence $M'$ is pseudo-null, thus
$q(M')=0.$  This proves (i).\\
 For any $M\in\mod(p),$ the
canonical decomposition
 $$0\subseteq {_p}M\subseteq {_{p^2}}M\subseteq\cdots \subseteq {_{p^m}}M=M$$
for some $m,$ induces a decomposition
 $$0\subseteq q({_p}M)\subseteq q({_{p^2}}M)\subseteq\cdots\subseteq  q({_{p^m}}M)=q(M)$$
with $$q({_{p^{j+1}}}M)/q({_{p^j}}M)\cong
q({_{p^{j+1}}}M/{_{p^j}}M)\cong q(\La/p)^{d_j},$$ where
$d_j=\rk_{\La/p}({_{p^{j+1}}}M/{_{p^j}}M)$ by the previous lemma.
Since this filtration can be refined easily to a decomposition
series of the desired kind, we are done.
\end{proof}

\begin{cor}\label{additiv}
The invariant $\mu$ is additive on short exact sequences of
\La-torsion modules.
\end{cor}

\begin{proof}
Since $\mu$ is additive on short exact sequences of $p$-torsion
modules by the proposition it suffices to reduce the general
statement to this case. Let
 $$\xymatrix@1{
    {\ 0 \ar[r] } &  {\ X \ar[r] } &  {\ Y \ar[r] } &  {\ Z \ar[r] } &  {\
    0 } \
 }$$ be a short exact sequence of \La-torsion modules. Choosing a
 number  $n$ such that the $p$-torsion parts of $X,Y$ and $Z$ are
 annihilated by $p^n,$ we obtain an exact sequence
  $$\xymatrix@1{
   {\ 0 \ar[r] } &  {\ \tor_\zp X \ar[r] } &  {\ \tor_\zp Y \ar[r] } &  {\ \tor_\zp Z \ar[r] } &  {\
   X/p^n\ar[r]^\varphi } & Y/p^n.
 }$$
Considering the exact, commutative diagram
 $$\xymatrix{
  0\ar[r]&{\tor_\zp X}\ar[r]\ar@{^(->}[d] & X/p^n\ar[r]\ar[d] & (X/\tor_\zp)/p^n\ar[r]\ar[d] & 0 \\
   0\ar[r]& {\tor_\zp Y}\ar[r] & Y/p^n\ar[r] & (Y/\tor_\zp)/p^n\ar[r] &
   0,
 }$$ we see that $\ker(\varphi)$ is pseudo-null by the following
 lemma, i.e.\
 $$\xymatrix@1{
   {\ 0 \ar[r] } &  {\ \tor_\zp X \ar[r] } &  {\ \tor_\zp Y \ar[r] } &  {\ \tor_\zp Z \ar[r] } &
   0
}$$ is exact $\mbox{mod }{\mathcal{PN}}.$
\end{proof}

\begin{lem}
\label{mu-vanish} Assume that $G$ is  a $p$-adic analytic group
without $p$-torsion such that both $\La=\La(G)$ and $\La/p$ are
integral. Let $M$ be a (not necessarily torsion) \La-module. Then
the following holds
 $$\mu(M/p)=0\:\:\:\Rightarrow\:\:\:\mu({_pM})=0.$$
 \end{lem}

\begin{proof}
Since $(\tor_\zp M)/p\subseteq M/p$ by the snake lemma, it
suffices to deal with the case that $M$ is \La-torsion. But then
the additivity of the $\mu$-invariant shows immediately that
$\mu({_pM})=\mu(M/p).$
\end{proof}

\begin{lem}
Assume that $G$ is  a $p$-adic analytic group without $p$-torsion
such that both $\La=\La(G)$ and $\La/p$ are integral. Let $M$ be a
\La-torsion module with $\tor_\zp M=0.$ Then
\begin{enumerate}
\item for any integer $n\geq1,$ the module $M/p^n$ is pseudo-null.
\item $\tor_\zp\E^1(M)=0.$
\end{enumerate}
\end{lem}

We will denote the annihilator in \La\ of an element $m\in M$ by
$ann_\La(m):=\{\lambda\in\La|\lambda m=0\}.$

\begin{proof}
Since there is a surjection
 $$\bigoplus \La/ann_\La(m_i)\twoheadrightarrow M$$
for a finite set of generators $m_i$ of $M,$ it suffices to prove
(i) in the case $M:=\La/I,$ where $I$ is a non-zero left ideal of
\La. As $M/p^n$ is \La-torsion we are done once we have shown the
vanishing of $\E^1_\La(M/p^n).$ But
 \begin{eqnarray*}
 \E^1_\La(M/p^n)&\cong&\E^0_{\La/p^n}(M/p^n)\\
                &\cong&\Hom_{\La/p^n}(M,\La/p^n)\\
                &\cong&\Hom_\La(\La/I,\La/p^n)=0.
 \end{eqnarray*}
Indeed, the vanishing of the latter module can be seen as follows:
let $\varphi:\La\to \La/p^n$ be a non-trivial homomorphism of
\La-modules which factors through $\La/I,$ i.e.\ $I\subseteq
ann_\La(x\mbox{ mod}p^n)$ with $x\equiv\varphi(1)\mbox{
mod}p^n.$\\

{\em Claim:} $ann_\La(x\mbox{ mod}p^n)\subseteq p\La.$\\

Let $\lambda\in ann_\La(x\mbox{ mod}p^n),$ i.e.\ $\lambda x=p^ny$
for some $y\in\La$ and let $n_o$ be the maximal integer with $x\in
p^{n_0}\La,$ i.e.\ $x=p^{n_0}x_0$ for some $x_0\in\La\setminus
p\La$ and $n_0<n.$ Since the multiplication by $p^{n_0}$ is
injective we obtain $\lambda x_o=p^{n-n_0}y\equiv 0\mbox{ mod}p.$
Hence $\lambda\in p\Lambda$ because $\Lambda/p$ is integral. This
proves the claim.\\
The fact that $_pM=0,$ implies $I\cap
p\Lambda=pI$ and regarding the claim it holds
 $$I=I\cap p\Lambda=pI=\ldots=p^mI$$
for any $m\geq0.$ Since $p^m$ tends to zero if $m$ goes to
infinity the ideal $I$ must be zero, a contradiction.\\
 The second statement results from the first one regarding the exact sequence
 $$\xymatrix@1{
    {\ 0=\E^1(M/p) \ar[r] } &  {\ \E^1(M) \ar[r]^{p} } &  {\ \E^1(M).  }  \
 }$$
\end{proof}

We finish this section with a ``structure theorem for the
$p$-torsion part of \La-modules."

\begin{thm}\label{p-structure}
Assume that $G$ is  a $p$-adic analytic group without $p$-torsion
such that both $\La=\La(G)$ and $\La/p$ are integral. Let $M$ be
in $\mod(p).$ Then there exist uniquely determined natural numbers
$n_1,\ldots,n_r$ and an isomorphism in $\mod(p)/{\mathcal{PN}}(p)$
 $$M\equiv \bigoplus_{1\leq i\leq r} \La/p^{n_i}\ \mbox{\rm  mod
 }{\mathcal{ PN}}(p).$$
\end{thm}

\begin{proof}
Let $m$ be minimal with the property: $_{p^{m+1}}M=M.$ The theorem
is proved using induction with respect to $m.$ The case $m=0$ is
just lemma \ref{fp-lemma}, so let $m$ be arbitrary. Again by lemma
\ref{fp-lemma} we are in the following situation:
 $$\xymatrix{
      &   &   & (\La/p)^d\ar@{^(->}[d]^\iota &   \\
   0\ar[r] & {_{p^m}}M\ar[r] & M\ar[r] & M/{_{p^m}}M\ar[r] & 0, \
 }$$
where $d$ is the $\La/p$-rank of $M/{_{p^m}}M$ and the cokernel of
$\iota$ is  pseudo-null. Replacing $M$ by the pull-back with
$\iota,$ we may assume that $M/{_{p^m}}M\cong(\La/p)^d.$ Since
$(\La/p^{m+1})^d$ is free in the category of
$\La/p^{m+1}$-modules, we obtain easily the following exact and
commutative diagram
 $$\xymatrix{
   0\ar[r] & (\La/p^m)^d \ar[r]\ar[d]^\varphi &(\La/p^{m+1})^d \ar[r]\ar[d]^\psi& (\La/p)^d\ar[r]\ar@{=}[d] & 0 \\
   0 \ar[r]& {_{p^m}}M\ar[r]\ar@{->>}[d] & M\ar[r]^{pr}\ar@{->>}[d] & (\La/p)^d\ar[r] & 0 \\
     & N\ar@{=}[r] & N &   &   \
 }$$
where $N$ is by definition the cokernel of $\psi$ respectively
$\varphi.$ First we will show that $\psi$ and hence also $\varphi$
is injective. Since $(\La/p^{m+1})^d$ -  being of projective
\La-dimension $1$ - does not contain any proper pseudo-null
\La-submodule, it suffices to prove that $\ker(\psi)$ is
pseudo-null. Assuming the contrary, i.e.\ that
$\mu(\ker(\psi))\neq0,$ it follows that
$\mu({_{p^{m+1}}}K/{_{p^{m}}}K)<d$ for the image $K$ of $\psi$
because for an arbitrary $p$-torsion \La-module $N$
$\rk_{\La/p}({_{p^{i+1}}}N/{_{p^{i}}}N)\geq\rk_{\La/p}({_{p^{i+2}}}N/{_{p^{i+1}}}N)$
holds for any $i\geq 0.$ But this contradicts the surjectivity of
$pr\circ\psi.$\\ To prove the theorem we only have to show that
$\varphi$ has a co-section in \linebreak
$\mod(p)/{\mathcal{PN}}(p),$ i.e.\ that the short exact sequence
in the left column splits. Indeed, then a section
$N\hookrightarrow{_{p^m}}M$ would give rise to a section
$N\hookrightarrow M,$ i.e.\ $M\cong N\oplus (\La/p^{m+1})^d,$ and
by the assumption of the induction $N$ is already of the desired
form. {\em  Here and in what follows we are arguing in the
quotient category $\mod(p)/{\mathcal{PN}}(p),$ though
  we omit the functor $q$ in the notation of maps and objects for simplicity.}\\
Again by this assumption, the module ${_{p^m}}M$ is isomorphic to
a module of the form $(\La/p^m)^{d'}\oplus\bigoplus_i
\La/p^{n_i},$  where $n_i<m.$ Assume first that $d=1.$  We claim
that the image of $\varphi$ is mapped surjectively onto one of the
factors $\La/p^m$ under the correspondent projection. Indeed, it
is easy to see that otherwise the image would be contained in
${_{p^{m-1}}}M,$ which contradicts the injectivity of $\varphi.$
Counting $\mu$-invariants, we see that $\varphi$ followed by the
 projection onto such a factor gives an isomorphism  and therefore induces the desired co-section. If $d>1$ we
make  the same procedure iteratively for every factor of
$(\La/p^m)^d$ after first splitting up the image of the previous
factor(s). The theorem follows because the uniqueness can be
deduced easily from the decomposition
  $$0\subseteq {_p}M\subseteq {_{p^2}}M\subseteq\cdots\subseteq  {_{p^m}}M=M$$
 counting $\La/p$-ranks.
\end{proof}

\begin{rem}
If one replaces $p\in\La(G)$ by any element $f$ in the center of
$\La(G)$ such that $(f):=\La(G)f$ is a prime ideal, i.e.\ such
that $\La(G)/(f)$ is integral, then one gets analogous results
concerning the ``$f$-torsion part" $\{m\in M| f^nm=0 \mbox{ for
some } n\}$ of $M.$ In particular, one obtains further invariants
$\mu_f$ for these prime elements.
\end{rem}
\POP
\PUSH{localduality.tex}%

\section{Spectral sequences for Iwasawa adjoints}

In the previous section we have seen that the functors $E^i(-)$
play an essential role in the dimension theory of \La-modules. In
this section we will mention some techniques  which  sometimes
allow to calculate these adjoints in applications when $\La$ is
the Iwasawa algebra $\La(G)$ of a profinite group $G.$ A part of
the results stems from  U. Jannsen (\cite{ja-is})  whom I would
like to thank heartily  for giving me his manuscript on a spectral
sequence for Iwasawa adjoints (\cite{ja-sps}).

 We shall write $\d(G)$ and $\C(G)$ for the categories of discrete and compact
\La-modules, respectively, whereas we denote the full
subcategories of cofinitely and finitely generated modules by
$\d_{cfg}(G)$ and $\C_{fg}(G)$, respectively.

Now, let $G=H\times \Gamma$ be the product of  profinite groups
$H$ and $\Gamma$. Assume that $\La(G)$ is Noetherian and that
$\Gamma$ is separable, i.e.\ it possesses a countable ordered
system of open normal subgroups $\Gamma_n$ as a  basis of open
neighborhoods of $1\in\Gamma$. Let $(\d_{cfg}(G))^\mathbb{N}$ be
the category of inverse systems in $\d_{cfg}(G)$ and consider the
left exact functor
 $$ T_\Gamma :  \d_{cfg}(G)\rightarrow (\d_{cfg}(G))^\mathbb{N},$$
 which sends $B$ to the inverse system
 $\{ B^{\Gamma_{n+1}}\stackrel{N_{\Gamma_n/\Gamma_{n+1}}}{\longrightarrow} B^{\Gamma_n}\},$
 and
 $$ \projlim{n}
 \Hom_{\La(H)}(-^\vee,\La(H)):(\d_{cfg}(G))^\mathbb{N}\rightarrow\mbox{$\La(G)$-Mod.}$$
Here the action of $\Gamma$ on $f\in\Hom_{\La(H)}(M,\La(H))$ for
$M\in\C(G)$ is defined by $(\gamma f)(m):=f(\gamma^{-1}m)$,
whereas $h\in H$ acts by the rule $(h f)(M):=f(m)h^{-1}$ as usual.

Since the category $(\d_{cfg}(G))^\mathbb{N}$ has enough
injectives, because $\d_{cfg}(G)$ has (\cite{ja-con}, Prop.\ 1.1),
we can form the right derived functors
 $$R^iT_\Gamma(B)=\{\H^i(\Gamma_{n+1},B)\stackrel{cor}{\longrightarrow}\H^i(\Gamma_{n},B)\}$$
 and $$R^i( \projlim{n} \Hom_{\La(H)}(B^\vee,\La(H))),$$
which  equals
 $$\projlim{n}R^i\Hom_{\La(H)}(B^\vee,\La(H))  $$
(cf.\ \cite{ja-con} Prop.\ 1.2, 1.3), if we restrict  ourselves to
elements of the subcategory $(\d')^\mathbb{N}$ where  $\d'$ is the
abelian  subcategory of $\d_{cfg}(G)$ consisting of
$\La(G)$-modules, which are cofinitely generated over $\La(H)$.
Indeed, in this case, the modules $\Hom_{\La(H)}(B_n^\vee,\La(H))$
are compact, i.e.\ the inverse limit functor is exact on the
corresponding inverse systems. Since $R^i\Hom_{\La(H)}(-,\La(H))$
extends the functors $\E^i_{\La(H)}(-)$ naturally from $\C(H)$ to
$\C_{fg}(G),$   we will  write also $\E^i_{\La(H)}(-)$ for the
first functor. Note that it  is endowed with a natural
$\Gamma$-action.

\begin{lem}\label{acylic}
The functor $T_\Gamma$ sends injectives to
$\projlim{n}\Hom_{\La(H)}(-^\vee,\La(H))$-acylics.
\end{lem}

\begin{proof}
It suffices to prove that $\zp\kl H\kr [\Gamma/\Gamma_n]$ is
$\Hom_{\La(H)}(-,\La(H))$-acyclic. But, for any resolution of
$\zp\kl H\kr [\Gamma/\Gamma_n]$ by $\La(G)$-projectives
 $$P^\bullet\rightarrow\zp\kl H\kr [\Gamma/\Gamma_n],$$ the sequence
 $$0\rightarrow\Hom_{\La(H)}(\zp\kl H\kr [\Gamma/\Gamma_n],\La(H))\rightarrow\Hom_{\La(H)}(P^\bullet,\La(H))$$
 is exact, because both, $\zp\kl H\kr [\Gamma/\Gamma_n]$ and the $P^i,$
 are projectives in $\C(H)$ (cf.\ \cite{nsw} (5.3.13)). The result follows by taking
 homology.
\end{proof}

 The Grothendieck spectral sequence for the
composition of the above functors gives

\begin{thm}

With notation as above,  there is a convergent cohomological
spectral sequence

$$ \projlim{n}
\E^i_{\La(H)}(\H^j(\Gamma_n,B)^\vee)\Rightarrow\E^{i+j}_{\La(G)}(B^\vee)$$
for any $B$ in $\d_{cfg}(G).$
\end{thm}

Note that all modules that occur in the spectral sequence are
compact $\La(G)$-modules.

\begin{proof}

The functor $\E^0_{\La(G)}(-)$ is the composition of the functors
$T_\Gamma$ and \linebreak $ \projlim{n}
\Hom_{\La(H)}(-^\vee,\La(H))$, because by  lemma \ref{com-lem} we
have isomorphisms of $\La(G)$-modules

 \begin{eqnarray*}
 \E^0_{\La(G)}(B^\vee)&=& \Hom_{\La(G)}(B^\vee,\zp\kl H\kr )\kl \Gamma\kr )\\
                 &=&\projlim{n}\Hom_{\zp\kl H\kr [\Gamma/\Gamma_n]}((B^\vee)_{\Gamma_n},\zp\kl H\kr [\Gamma/\Gamma_n])\\
                 &=&\projlim{n}\Hom_{\La(H)}((B^{\Gamma_n})^\vee,\La(H)).
 \end{eqnarray*}

Now the result follows by lemma \ref{acylic}.

\end{proof}

Recall that there is a canonical $\La(H)$-homomorphism
$$\pi_n:\zp\kl H\kr [\Gamma/\Gamma_n]\rightarrow\zp\kl H\kr ,
\sum_{g\in \Gamma/\Gamma_n} a_g\ g\Gamma_n\mapsto a_1,$$ and, for
any $m\geq n,$ a canonical $\La(G)$-homomorphism $p_{m,n}:\zp\kl
H\kr [\Gamma/\Gamma_m]\rightarrow\zp\kl H\kr [\Gamma/\Gamma_n]$
which sums up the coefficients of the same $\Gamma_n$-cosets.

\begin{lem}\label{com-lem}
The homomorphisms $\pi_n$ and $p_{m,n}$ induce a
commutative\linebreak diagram of $\La(G)$-modules:
 $$\xymatrix{
  { \Hom_{\zp\kl H\kr [\Gamma/\Gamma_m]}(M_{\Gamma_m},\zp\kl H\kr [\Gamma/\Gamma_m] )}\ar[r]^(.5){(\pi_m)_\ast}_(.5)\simeq \ar[d]^{(p_{m,n})_\ast}& {\phantom{a}\Hom_{\zp\kl H\kr }(M_{\Gamma_m},\zp\kl H\kr )}\ar[d]^{N_{\Gamma_n/\Gamma_m}} \\
  {\Hom_{\zp\kl H\kr [\Gamma/\Gamma_n]}(M_{\Gamma_n},\zp\kl H\kr [\Gamma/\Gamma_n] )}\ar[r]^(.5){(\pi_n)_\ast}_(.5)\simeq &{\phantom{a}\Hom_{\zp\kl H\kr }(M_{\Gamma_n},\zp\kl H\kr )} \
 }$$
\end{lem}

\begin{proof}
It is easily checked that the diagram commutes and that the
inverse of $(\pi_n)_\ast$ is given by $\psi\mapsto
(m\mapsto\sum_{g\in\Gamma/\Gamma_n}\psi(g^{-1}m) g\Gamma_n)$.
(Note that the $\Gamma$-invariance of a homomorphism $\phi(m)=\sum
\phi(m)_g g\Gamma_n$ is equivalent to
$\phi(\gamma^{-1}m)_1=\phi(m)_\gamma$ for all $\gamma\in\Gamma$.)
Recalling that $\gamma\in\Gamma$ acts by  $(\gamma
\phi)(m):=\phi(\gamma^{-1}m)$ on $\Hom_{\zp\kl H\kr
[\Gamma/\Gamma_n]}(M_{\Gamma_n},\zp\kl H\kr [\Gamma/\Gamma_n] ),$
it is also immediate that $(\pi_n)_\ast$ is $\La(G)$-invariant.
\end{proof}

\begin{cor}\label{zp-factor}
If $\Gamma$ contains an open subgroup of index prime to $p$ and
isomorphic to $\zp,$  then there is a long exact sequence of
$\La(G)$-modules

{\footnotesize $$ \xymatrix@1@C=12pt{
   { \projlimsc{n}^{\phantom{\stackrel{ M}{N}}}\hspace{-6pt}\E_{\La(H)}^i(M_{\Gamma_n})} \ar[r]  &
    {\ \E^i_{\La(G)}(M)  \ar[r] } &
     {\ \projlimsc{n}^{\phantom{\stackrel{ M}{N}}}\hspace{-6pt}\E^{i-1}_{\La(H)}(M^{\Gamma_n}) \ar[r] } &
      {\ \projlimsc{n}^{\phantom{\stackrel{ M}{N}}}\hspace{-6pt}\E_{\La(H)}^{i+1}(M_{\Gamma_n}) \ar[r] } &
   {\   \E^{i+1}_{\La(G)}(M) }
   }
$$}
\end{cor}

Now we are going to present  further spectral sequences due to U.
Jannsen which were in some sense the models for the first one
proved in this section. The next one describes the Iwasawa
adjoints of certain cohomology groups associated with $p$-adic
representations. So let $G$ be a profinite group and $G_{\infty}$
a closed subgroup, such that its quotient has a countable basis of
neighbourhoods of identity, i.e.\ there is a countable family
$G_n,$ $G_{\infty}\subseteq G_n\subseteq G,$ with $\bigcap_n
G_n=G_{\infty}.$ Furthermore, let $A=(\Qp/\zp)^r$ for some $r\geq
1$ with some continuous action of $G$. We shall write
 $$T_pA=\Hom(\qp/\zp,A)\cong\projlim{m} {_{p^m}A}$$
for the Tate module of $A.$ Then there is the following convergent
spectral sequence (\cite{ja-sps}):

\begin{thm}(Jannsen)\label{ja-ss}
$$E_2^{p,q}=\E^p(\H^q(G_{\infty},A)^\vee)\Rightarrow\projlim{n}
\H^{p+q}(G_n,T_p A)$$
\end{thm}

\begin{cor}
Assume $\cd_p(G)\leq 2$. Then the exact sequence of low degrees
degenerates to
  {$$\xymatrix@R-2.5pc{
  0\ar[r] & {\E^1(A(k_\infty)^\vee)\ar[r]} & {\projlim{n}^{\phantom{\stackrel{ M}{N}}}\hspace{-6pt}\H^1(G_n,T_p A)\ar[r]} & {\E^0(\H^1(G_{\infty},A)^\vee)}
  \ar[r]& {\phantom{0.}} \\
& {\E^2(A(k_\infty)^\vee)} \ar[r]&
   {\ker(\projlim{n}^{\phantom{\stackrel{\stackrel{ m}{N}}{N}}}\hspace{-6pt}\H^2(G_n,T_p
   A)}\ar[r]&{\E^0(\H^2(G_{\infty},A)^\vee))\ar[r]}&{\phantom{0.}} \\
 & & {\E^1(\H^1(G_{\infty},A)^\vee)\ar[r]} & {\E^3(A(k_\infty)^\vee)\ar[r]} &
  0.}$$}
\end{cor}

The next result, which relates the (compact) $\Lambda$-modules
$\E^i(M)$ to the discrete $G$-modules
 $$D_i(M^\vee):=\dirlim{U\subseteq_o G} \H^i(U,M^\vee)^\ast\ ,\ i\geq 0,$$ is derived by some spectral
sequences, too, but we only state the associated long,
respectively short, exact sequences:

\begin{thm}\label{ja-tor} (Jannsen) Let $G$ be a profinite group such that
$\Lambda(G)$ is Noetherian. Then, for any finitely generated
$\Lambda$-module $M,$ there are functorial exact sequences
\begin{enumerate}
\item[\rm(i)]

{ $$\xymatrix@1{
   {\ 0 \ar[r] } &  {\ D_i(M^\vee)\otimes_\zp \Qp/\zp \ar[r] } &  {\ \E^i(M)^\vee \ar[r] } &  {\ \mathrm{tor}_\zp D_{i-1}(M^\vee) \ar[r] } &  {\
   0, }
}$$} for all $i,$ where by definition $D_{i}(M^\vee)=0$ for $i<0.$
\item[\rm(ii)]

{ $$\xymatrix@1{
  { \ar[r] } &{\E^i(M)^\vee \ar[r] } &  {\dirlim{m}^{\phantom{\stackrel{ M}{N}}}\hspace{-6pt} D_i(_{p^m}(M^\vee)) \ar[r] } &  {\dirlim{m}^{\phantom{\stackrel{ M}{N}}}\hspace{-6pt} D_{i-2}(M^\vee/p^m) \ar[r] } &  {\E^{i-1}(M)^\vee
   }\ar[r]& ,
}$$}
\end{enumerate}
and the following isomorphisms
\begin{enumerate}
\item[\rm (iii)] $\E^i(M/\mathrm{tor}_\zp M)^\vee\cong \dirlim{m}
D_i(_{p^m}(M^\vee)),$
\item[\rm (iv)] $\E^i(\mathrm{tor}_\zp M)^\vee\cong \dirlim{m}
D_{i-1}(M^\vee/p^m).$
\end{enumerate}
\end{thm}

\begin{proof}
See \cite{ja-is} 2.1, 2.2 or \cite{nsw} theorem 5.4.12.
\end{proof}

For the duality theorem in the next section the following
corollary will be crucial.

\begin{cor}\label{ja-cor}
Assume that $G$ is a duality group at $p$ of dimension $n$ with
dualizing module $D_n^{(p)}=\dirlim{m}
D_n(\mathbf{Z}/p^m\mathbf{Z}).$ Then the following holds:
\begin{enumerate}
\item If $M$ is $\Lambda$-module which is free of finite rank as
$\zp$-module, then
 $$\E^i(M)^\vee\cong\left\{\begin{array}{cl}
   \dirlim{m} D_n((M/p^m)^\vee)\cong M\otimes_\zp D_n^{(p)} & \mbox{if  } i=n, \\
   0  & \mbox{otherwise.}
 \end{array}\right.$$
\item If $N$ is a finite $p$-primary $\Lambda$-module, then
$$ \E^i(N)^\vee\cong\left\{\begin{array}{cc}
  \Hom_\zp (N^\vee,D_n^{(p)}) & \mbox{ if  } i=n+1, \\
  0 &\mbox{otherwise.}
\end{array}\right.$$
\end{enumerate}
\end{cor}

\begin{proof}
See \cite{ja-is} 2.6 or \cite{nsw} 5.4.14.
\end{proof}

\begin{prop}\label{ind-dimen}
Let $G$ be a compact $p$-adic analytic group without $p$-torsion,
$H\subseteq G$  a closed subgroup and $M$ a finitely generated
$\La(H)$-module. If $d_{\La(G)}$ (resp.\ $d_{\La(H)}$) denotes the
(projective or $\delta$-) dimension of $\La(G)$ (resp.\ $\La(H)$),
then the following holds:
\begin{enumerate}
\item $j_{\La(G)}(\Ind^H_G M)=j_{\La(H)}(M),$
\item $\delta_{\La(G)}(\Ind^H_G M)=\delta_{\La(H)}(M) + d_{\La(G)}-d_{\La(H)},$
\item $\pd_{\La(G)}(\Ind^H_G M)=\pd_{\La(H)}(M).$
\end{enumerate}
\end{prop}

\begin{proof}
This is a consequence of \ref{inducedExt}, \ref{basic1},
\ref{dim-formel}, and \ref{pdviaExt}.
\end{proof}

\begin{lem}\label{corlemma}
Assume that $G=H\times\Gamma$ is a $p$-adic Lie group without
$p$-torsion where $\Gamma$ contains an open subgroup of index
prime to $p$ which is isomorphic to $\zp.$  Let $M\in\C(G)$ be
finitely generated and torsion as $\La(H)$-module. Then $M$ is a
pseudo-null $\La(G)$-module.
\end{lem}

\begin{proof}
By the corollary \ref{zp-factor}, there is an exact sequence
 $$
 \xymatrix@1{
    {\ 0 \ar[r] } &  {\ \projlim{n} \E^1_{\La(H)}(M_{\Gamma_n} )} \ar[r]  &{\ \E^1_{\La(G)}(M)} \ar[r]
    & {\ \projlim{n}\E^0_{\La(H)}(M^{\Gamma_n})=0. }
 }
$$ So, if we can show that the left term vanishes, we are done,
because then $\E^1\E^1(M)=0=\E^0\E^0(M).$ Consider the commutative
exact diagram
 $$\xymatrix{
  M\ar[r]^{\omega_n}\ar@{=}[d] & M\ar[r]\ar[d]^{\frac{\omega_m}{\omega_n}} & M_{\Gamma_n}\ar[r]\ar[d]^{\frac{\omega_m}{\omega_n}} & 0 \\
   M \ar[r]^{\omega_m}  & M\ar[r]  & M_{\Gamma_m}\ar[r] & 0,
 }
 $$
 where $\omega_n=\gamma^{p^n}-1$ for some generator $\gamma$ of
 $\zp\subseteq \Gamma$. Since $M$ is assumed to be $\La(H)$-torsion, we get the
 commutative exact diagram
  $$\xymatrix{
    0\ar[r] & {\E^1_{\La(H)}(M_{\Gamma_m})}\ar[r]\ar[d]^{\frac{\omega_m}{\omega_n}} & {\E^1_{\La(H)}(M)\ar[r]^{\omega_m}\ar[d]^{\frac{\omega_m}{\omega_n}} }&{\E^1_{\La(H)}(M)}\ar@{=}[d] \\
  0\ar[r]& {\E^1_{\La(H)}(M_{\Gamma_n})}\ar[r] & {\E^1_{\La(H)}(M)\ar[r]^{\omega_n} } & {\E^1_{\La(H)}(M).} \\
  }$$
Passing to the limit, we obtain $$\projlim{n}
\E^1_{\La(H)}(M_{\Gamma_n}
)\subseteq\projlim{n}\E^1_{\La(H)}(M)=\projlim{n}\bigcap_{m\geq
n}\frac{\omega_m}{\omega_n}\E^1_{\La(H)}(M)=0,$$ because
$\frac{\omega_m}{\omega_n}$ tends to zero.

\end{proof}

\begin{rem}\label{HdimvergleichG}
The same arguments show that $\delta_G(M)\leq\delta_H(M)$ for any
finitely generated $\La(H)$-module $M.$
\end{rem}

Besides  the case  $G=\mathbb{Z}_p^d$ these results apply also to
the following situation where $G$ is an open subgroup of
$Gl_d(\zp),$  $d$ is prime to $p,$  such that the determinant
takes values in
$\Gamma:=\det(G)\subseteq\zp\subseteq\mathbb{Z}_p^\ast.$ at least
if . Indeed, we have the following exact commutative diagram
  $$\xymatrix{
    1\ar[r] & Sl_d(\zp)\ar[r] & Gl_d(\zp)\ar[r]^{\det} & {\mathbb{Z}_p^\ast}\ar[r]& 1 \\
    1\ar[r] & Sl_d(\zp)\cap G\ar[r]\ar@{^{(}->}[u] & G\ar@{^{(}->}[u]\ar[r] & {\Gamma}\ar@{^{(}->}[u]\ar[r] & 1,
  }
$$ in which the lower sequence possesses  the following splitting
$$s:\Gamma\cong\zp\rightarrow G,\
a\mapsto\left(\begin{array}{cccc}
  a^{1/d} & &   &   \\
    & a^{1/d} &   &   \\
    &   & \ddots &   \\
    &   &    & a^{1/d}
\end{array}\right)$$
(Note that $\Gamma\cong\zp$ is considered as subgroup of the units
and that $\zp$ is uniquely $d$-divisible. Furthermore, if the
image of this homomorphism is not contained in $G,$ we just apply
the theory to an open subgroup $U$ of $G$   which fulfills this
condition with respect to $\det(U)$ and contains $H:=Sl_d(\zp)\cap
G.$ Such $U$ always exists because $Gl_d(\zp)$ is $p$-adic
analytic, i.e.\ the lower $p$-series forms a basic of
neighborhoods of the neutral element. Hence at least for some $m$
the image of $p^m\Gamma$ is contained in $G:$
$s(a^{p^m})=s(a)^{p^m}\in P_m(Gl_d)\subseteq G.$ Take
$U:=\det^{-1}(p^m\Gamma)\cap G.$) Since the splitting takes values
in the center of $G$, we get a presentation of $G$ as the direct
product $G=H\times \zp$.

\section{Local Duality}
\label{localdualitysec}

In this and the following section let $\La=\La(G)=\zp\kl G\kr $ be
the completed group algebra over $\zp$, {\em where $G$ is a
pro-$p$ Poincar\'{e} group, such that \La\ is Noetherian,} $\M$ the
maximal ideal of \La\ and $k=\La/\M\cong\mathbb{F}_p$ its finite
residue class field. It is well known that the global homological
dimension of \La\ is $d=\cd(G)+1$. By \Mod we denote the category
of (abstract) modules over the (abstract) ring \La\ and we write
\mod for the full subcategory of finitely generated modules. In
the sequel we will use frequently the equivalence of the latter
category with the category of finitely generated compact  modules.

\begin{defn}
For a finitely generated \La-module $M$, we define the depth  by
 $$ \depth(M):=\min\{i\mid \Ext_{\La}^i(k,M)\neq0\}.$$
\end{defn}

Recall that for a commutative Noetherian ring \La\ the $I$-depth
$depth_I(M)$   of a  finitely generated \La-module $M$ with
respect to an ideal $I$ is the maximal length of a $M$-regular
sequence in $I.$ For a local ring the $\depth(M)$ is
 $depth_\M (M),$ while the grade defined in
\ref{defAuslander} is  $j(M)=depth_{ann(M)}(\La),$ where $ann(M)$
is the annihilator of $M$ in \La.

We consider the additive functor
$\Gamma_{\M}(-):\Mod\rightarrow\Mod$ defined by
$\Gamma_{\M}(M):=\{x\in M\mid \M^l x=0\mbox{ for some l } \}$ and
state some basic properties:

\begin{lem}
\begin{enumerate}
\item $\Gamma_{\M}(M)=\dirlim{l} \Hom_{\La}(\La/\M^l ,M),$\\
in particular, the functor $\Gamma_{\M}(-)$ is left exact.
\item The restriction of \/ $\Gamma_{\M}$ to \mod equals $\t0$, i.e.\
$\Gamma_{\M}(M)$ is the maximal finite submodule of $M$, if the
latter module is finitely generated.
\end{enumerate}
\end{lem}

\begin{proof}
Since $\Hom_{\La}(\La/\M^l,M)=\{x\in M\mid \M^lx=0\},$ the first
statement is obvious. If $M$ is finitely generated, there is some
$l$ such that $\M^l\t0(M)=0$, i.e.\
$\t0(M)\subseteq\Gamma_{\M}(M).$ On the other hand $\La/\M^l$ is a
finite ring. Therefore $\La x\subseteq\t0(M)$ holds for any $x\in
\Gamma_{\M}(M).$
\end{proof}

Since \Mod has sufficiently many injectives, we can form the right
derived functors
  $$\H_{\M}^i(-)=R^i\Gamma_{\M}(-)=\dirlim{l}\Ext^i_{\La}(\La/\M^l,M)$$
  (noting the exactness of direct limits in \Mod).
We write $$\Mod_{\M},\ \mod_{\M}$$ for the full subcategory of
\Mod, \mod respectively, consisting of those modules $M,$ for
which $\H^0_{\M}(M)=M$ holds. $${\mathcal{ D}}(\Mod)\
(\mbox{resp.\ } {\mathcal{ C}}(\Mod))$$ means the category of
discrete (resp.\ compact) \La-modules, where \La\ is endowed with
its canonical $(\m,I)$-topology.

\begin{lem}
 $\H^i_{\M}(-)$ commutes with direct limits.
\end{lem}

\begin{proof}
Choose a resolution $P_{\bullet}$ of $\La/\M^l$ by finitely
generated projectives in order to calculate
$\Ext^i_{\La}(\La/\M^l,M).$ Since $\Hom_\La (P_j,-)$ commutes with
direct limits (as $P_j$ is finitely generated, i.e.\ any
homomorphism $\phi:P_j\rightarrow\dirlim{i} M_i$ factors over some
$M_i$), $\Ext^i_{\La}(\La/\M^l,-)$ does also and the lemma
follows.
\end{proof}

\begin{prop}
The forgetful functor defines an equivalence of categories
 $${\mathcal{ D}}(\Mod)\cong\Mod_{\M}.$$
\end{prop}

\begin{proof}
Both categories consists exactly of direct limits of finite
modules (cf.\ \cite[Prop.\ (5.2.4)]{nsw} for ${\mathcal{
D}}(\Mod)$).
\end{proof}

\begin{lem}\label{depth}
\begin{enumerate}
\item $\H^i_\M(\Mod)\subseteq\Mod_\M\cong{\mathcal{ D}}(\Mod)$ $\fa i\geq 0.$
\item For any $M\in\mod$, it holds $\depth(M)=\min\{i\mid\H^i_\M(M)\neq
0\}.$
\item \label{depth1} $\depth(\La)=d$ and $\H^d_\M(\La)=\La^\vee.$
\item $\Hom_\La (M,\H^d_\M(\La))\cong M^\vee \fa M$ in $\Mod_\M$ or in
\mod, in particular, $\H^d_\M(\La)$ is an injective \La-module.

\end{enumerate}
\end{lem}

\begin{proof}
Since $\H^i_\M(-)$ are the derived functors of $\H^0_\M(-),$ it
suffices to prove (i) for the latter functor. But in this case the
statement holds just by definition.

Now we will prove (ii) and set $k=\min\{i\mid\H^i_\M(M)\neq 0\}.$
Since \linebreak $\Ext^i_\La (\La/\M^l, M)=0$ for all
$i<\depth(M)$ (note that $\La/\M^l$ has a finite composition
series with subquotients isomorphic to $k$), it holds
$\depth(M)\leq k.$ So we only have to prove that $\H^j_\M(M)\neq
0$ for $j=\depth(M)<\infty.$ But the short exact sequences
 $$\xymatrix@1{
   {\ 0 \ar[r] } &  {\ \M/\M^l  \ar[r] } &  {\ \La/\M^l \ar[r] } &  {\ k \ar[r] } &  {\
   0 }
}$$ induce the long exact sequences
 $$\xymatrix@1{
   {\ 0=\Ext^{j-1}_\La(\M/\M^l,M) \ar[r] } &  {\ \Ext^j_\La(k,M) \ar[r] } &  {\ \Ext^j_\La(\La/\M^l,M) \ar[r] } &  {\ \cdots , }
}$$
i.e.\ $0\neq\Ext^j_\La(k,M)\subseteq\H^j_\M(M).$

Using \ref{ja-cor} and denoting the character of the dualizing
module by $\chi,$ we calculate
$$\H^i_\M(\La)=\dirlim{l}\E^i(\La/\M^l)=\left\{\begin{array}{cl}
  \dirlim{l} ({\La/\M^l}(\chi))^\vee=\La^\vee & \mbox{if } i=d \\
  0  & \mbox{otherwise,}
\end{array} \right.$$
whence (iii) follows. In order to prove (iv) first let $M$ be in
$\Mod_\M$, i.e.\ $M=\dirlim{i} M_i$ for some finite \La-modules
$M_i.$ Then, noting that $M_i$ is a $\La/\M^{l(i)}$-module for
some $l(i)$ and using the adjunction of ``$\Hom$ and $\otimes$",
 \begin{eqnarray*}
 \Hom_\La(M,\H^d_\M(\La))&=&\Hom_\La(\dirlim{i}M_i,\dirlim{l}
 ({\La/\M^l})^\vee)\\
 &=&\projlim{i}\Hom_\La(M_i,\dirlim{l}({\La/\M^l})^\vee)\\
 &=&\projlim{i}\Hom_\La(M_i,\Hom_\zp(\La/\M^{l(i)},\mathbb{Q}_p/\zp))\\
 &=&\projlim{i}\Hom_\zp(M_i,\mathbb{Q}_p/\zp)\\
 &=&M^\vee.
\end{eqnarray*}

Now let $M$ be in \mod. Then, noting that $\Hom_\La(M,-)$ commutes
with direct limits, because $M$ is finitely generated,
\begin{eqnarray*}
 \Hom_\La(M,\H^d_\M(\La))&=&\Hom_\La(M,\dirlim{l}
 {(\La/\M^l)}^\vee)\\
 &=&\dirlim{l}\Hom_\La(M,({\La/\M^l})^\vee)\\
 &=&\dirlim l\Hom_\La(M/\M^l,\Hom_\zp(\La/\M^{l},\mathbb{Q}_p/\zp))\\
 &=&\dirlim l\Hom_\zp(M/\M^l,\mathbb{Q}_p/\zp)\\
 &=&M^\vee.
\end{eqnarray*}
\end{proof}

After this technical preparations we are able to prove the
following

\begin{thm}\label{localduality}
Let $G$ be a pro-p Poincar\'{e} group  with $d:=\cd(G)+1<\infty$ and
such that $\La=\La(G)$ is Noetherian. Then, for any $M\in\mod,$

 $$\E^i(M)\cong\Hom_\La(\H^{d-i}_\M(M),\H^{d}_\M(\La))\cong\H^{d-i}_\M(M)^\vee=:T^i(M).$$
\end{thm}

\begin{proof}
Consider the right exact contravariant additive functor \linebreak
$T^0(-)=\H^{d}_\M(M)^\vee$ on \mod (note that $\H^{i}_\M(M)=0\fa
i>d$ as \La\ has global dimension $d$). By \cite[Thm. 3.36 and
Remarks ]{rot} there is a natural equivalence of functors
 $$T^0(-)\cong\Hom_\La(-,T^0(\La))\cong\Hom_\La(-,\La)$$
on \mod. Therefore, it suffices to show that the functors $T^i(-)$
are the left derived functors of $T^0(-)$. But $\{T^i(-)\}_{i\geq
0}$ is a universal $\delta$-functor because they are effaceable by
projectives in \mod (Since $T^0$ is additive, it is sufficient to
verify that $\H^i_\M(\La)=0$ for all $i<d,$ which is done by lemma
\ref{depth} (iii)).
\end{proof}
\POP
\PUSH{ab-equality.tex}%
\section{Auslander-Buchsbaum equality}

In this section we assume the same conditions on \La\ as in the
previous one and, under this conditions, we are going to prove the
Auslander-Buchsbaum equality
  $$\pd(M)+\depth(M)=\depth(\La)$$
for all $M\in\mod.$ In the theory of commutative local rings this
equality can be proved using {\em regular sequences}. Since this
concept is lacking in the non-commutative theory, we will have to
replace it by homological methods, i.e.\ we will work in derived
categories. Our proof is analogous to J\o rgensen's proof of the
Auslander-Buchsbaum equality in the case of (non-commutative)
graded algebras over a field (cf.\ \cite{joer}).

First, we recall the definitions of total Hom and total tensor
product. Let $X,Y\in\K(\Mod)$ and define

 $$(\Hom_\La(X,Y))^n=\prod_{i\in\mathbb{Z}}\Hom_\La(X^i,Y^{i+n}),\
 d^n=\prod_{i} (d_X^{i-1}+(-1)^{n+1}d_Y^{i+n})$$
and
 $$(X\otimes_\La Y)^n=\bigoplus_{i+j=n}X^i\otimes_\La Y^j,\
 d^n=\bigoplus_{i+j=n}(d_X^i\otimes 1 + (-1)^n\otimes d_Y^j).$$
 They become bifunctors

 $$\Hom_\La(-,-):\K(\Mod)^{op}\times\K(\Mod)\rightarrow\K(\zpMod),$$
 $$-\otimes_\La -:\K(\mbox{Mod-\La})\times\K(\Mod)\rightarrow\K(\zpMod),$$
where we denote by $\mbox{Mod-\La}$  the category of right
\La-modules. Note that the latter category is  equivalent to \Mod
due to the involution on the group algebra \La. Moreover, if $Y$
is a complex of bi-modules, then the values of $\Hom_\La(-,Y)$ are
in $\K(\mbox{Mod-\La})$, if $X$ is a complex of bi-modules, then
$X\otimes_\La -$ has values in $\K(\Mod)$.

Since \Mod has enough projectives, the derived functors exist
(cf.\ \cite[Chap.\ I, Theorem 5.1]{hart} or \cite[Thm
10.5.6]{weibel}):
 $$\rhom_\La(-,-):\D^-(\Mod)^{op}\times\D(\Mod)\rightarrow\D(\zpMod),$$
respectively
 $$\rhom_\La(-,-):\D^-(\Mod)^{op}\times\D(\mbox{\La-Mod-\La})\rightarrow\D(\mbox{Mod-\La})$$
and
 $$-\Lotimes_\La -:\D(\mbox{Mod-\La})\times\D^-(\Mod)\rightarrow\D(\zpMod),$$
respectively
 $$-\Lotimes_\La -:\D(\mbox{\La-Mod-\La})\times\D^-(\Mod)\rightarrow\D(\Mod).$$
$\rhom$, respectively $\Lotimes$, is computed via a projective
resolution in the first, respectively second variable.

\begin{prop}\label{homtensor}
Let $Y\in\D^b(\mbox{\La-Mod-\La}),$ $Z\in\D^b(\Mod)$ and let
\linebreak $X\in\D^b(\Mod)$ be a bounded complex which is
quasi-isomorphic to a bounded complex consisting of finitely
generated free \La-modules. Then
 $$\rhom_\La(X,Y\Lotimes_\La Z)\cong\rhom_\La(X,Y)\Lotimes_\La Z.$$
\end{prop}

\begin{proof}
(See \cite[Proposition 2.1]{joer} for the case of graded algebras
over a field.)

Replacing $X$ with a quasi-isomorphic complex $L\in\D^b(\Mod)$
consisting of finitely generated free \La-modules and replacing
$Z$ with a quasi-isomorphic complex $F\in\D^-(\Mod)$ consisting of
projectives, we see that  we have to prove
 $$\Hom_\La(L, Y\otimes_\La F)=\Hom_\La(L,Y)\otimes_\La F.$$
But due to the boundedness condition and the fact that $L$
consists of finitely generated free modules, the $n$th module on
either side becomes

 $$\bigoplus_{i,j} \Hom_\La(L^i,Y^j)\otimes_\La F^{n+i-j}$$ while
 the differentials on each summand $\Hom_\La(L^i,Y^j)\otimes_\La
 F^{n+i-j}$ are given by

 \begin{eqnarray*}
 &d_L^{i-1}\otimes 1 + d_Y^j\otimes (-1)^{j-i-1} + (-1)^{n}\otimes
 d_F^{n+i-j},& \ \mathrm{respectively}\\
 &d_L^{i-1}\otimes 1 + d_Y^j\otimes (-1)^{n} + (-1)^{i+1}\otimes
 d_F^{n+i-j}.&
 \end{eqnarray*}
We will construct an isomorphism between the two complexes: If the
minimal non-zero module of each of the complexes is
$\Hom(L^{i_0},Y^{j_0})\otimes\La F^{n_0+i_0-j_0}$, then the
multiplication by suitable signs on the summands associated to the
triple of indices $(a,b,c)=(i,j,n+i-j)$ defines an  isomorphism of
complexes. For example, we can choose these signs by the following
rules, which determine them uniquely:
\begin{enumerate}
\item $sign((i_0,j_0,n_0+i_0-j_0))=1,$
\item $sign((a+1,b,c))=sign(a,b,c),$
\item $sign((a,b+1,c))=(-1)^csign((a,b,c)),$
\item $sign((,a,b,c+1))=(-1)^{c+b+1}sign((a,b,c)).$
\end{enumerate}

\end{proof}

In the proof of the next theorem we use the notation

 \xymatrix{
  \sigma_{\geq n}(Y):=\cdots \ar[r] & 0 \ar[r] & Y^n/{\rm im}(Y^{n-1}) \ar[r] & Y^{n+1} \ar[r] & Y^{n+2} \ar[r] &
  \cdots}
 \noindent for the truncation of a complex $Y$ at the degree $n.$

\begin{thm} (Auslander-Buchsbaum equality) \label{ab-equa}
For any $M\in\mod$, it holds
 $$\pd_\La(M)+\depth_\La(M)=\depth_\La (\La).$$
\end{thm}

\begin{proof}(See \cite[Thm 3.2]{joer} for the case of graded algebras
over a field.)

 Regard $k,\ M,\ \La$ as complexes concentrated in
degree zero. Then the invariants in question are related to each
other by the following isomorphism
 $$\rhom_\La(k,M)\cong\rhom_\La(k,\La\Lotimes_\La
 M)\cong\rhom_\La(k,\La)\Lotimes_\La M,$$
where we use proposition \ref{homtensor}. Choosing a minimal free
resolution $L$ of  $M$ (see Appendix) and noting that the
truncation
 $$T=\sigma_{\geq d}(\rhom_\La(k,\La))$$
 is quasi-isomorphic to $\rhom_\La(k,\La),$ we can replace the
 right term by $T\otimes_\La L.$

  The lowest non-zero module in $T$
 is $T^d$ with $d=\depth(\La)$ while the lowest non-zero module in
 $L$ is $L^{-\pd(M)}$ according to Appendix, corollary \ref{mini-res-cor} . So the lowest
 non-zero module in $T\otimes_\La L$ becomes $(T\otimes_\La
 L)^{d-\pd(M)}=T^d\otimes_\La L^{-\pd(M)}.$ Obviously,
 $\depth(M)\geq d-\pd(M),$ because $\depth(M)=\min\{i\mid
 \H^i(\rhom_\La(k,M)\neq0\}.$ So we need to see that
 $\H^{d-\pd(M)}(T\otimes_\La L)$ is nonzero.

 However, $k\cong\Ext_\La^d(k,\La)=\ker(d_T^d)\subseteq T^d$ and
 the ``beginning" of the complex $T\otimes_\La L$ looks like

 \xymatrix@1{
    {\ 0 \ar[r] } &  {\ T^d\otimes_\La L^{-\pd(M)} \ar[r] } &  {\ T^d\otimes_\La L^{-\pd(M)+1}\oplus T^{d+1}\otimes_\La L^{-\pd(M)} \ar[r] } &  {\ \cdots  } . \
 }

Now it holds that
 $$0\neq\ker(d_T^d)\otimes_\La L^{-\pd(M)}\subseteq
 \ker(d_{T\otimes L}^{d-\pd(M)})=\H^{d-\pd(M)}(T\otimes_\La L).$$
Indeed, for $t\otimes l\in\ker(d_T^d)\otimes_\La L^{-\pd(M)},$ we
have $$d_{T\otimes L}^{d-\pd(M)}(t\otimes l)=d_T^d(t)\otimes
l+(-1)^{d-\pd(M)}t\otimes d_L^{-\pd(M)}(l).$$  The first summand
is zero because $t\in \ker(d_T^d)$ while, due to the minimality of
$L$ (cf.\ Appendix, proposition \ref{minimal res} (ii)), the
second one lies in $\ker(d_T^d)\otimes_\La \M
L^{d-\pd(M)+1}\cong\La/\M\otimes_\La \M L^{d-\pd(M)+1}=0. $
\end{proof}

\begin{cor}\label{pdviaExt}
If $M$ is a finitely generated \La-module, then
 $$\pd(M)=\max\{i\mid \E^i(M)\neq 0\}.$$
\end{cor}

\begin{proof}
Using lemma \ref{depth} (ii) and local duality, we get
\begin{eqnarray*}
\pd(M)&=& d-\depth(M)\\
      &=& d-\min\{i\mid\H^i_\M(M)\neq0\}\\
      &=& \max\{i\mid \E^i(M)\neq 0\}.\\
\end{eqnarray*}
\end{proof}

\begin{rem}
The statement of the last corollary holds over an arbitrary
Noetherian ring for a finitely generated modules $M$ with {\em
finite} projective dimension $\pd_\La M$  and can be proven
directly in the following way. Consider a projective resolution of
minimal length
 $$\xymatrix@1{
    {\ 0 \ar[r] } &  {\ P_n \ar[r]^{d_n} } &  {\ P_{n-1} \ar[r]^{d_{n-1}} } &  {\ \cdots \ar[r] } &  {\
    P_0\ar[r] } & M\ar[r]& 0. \
 }$$
Then the $(n-1)$th syzygy $K=\ker(d_{n-2})$ has projective
dimension $\pd_\La K=1,$ i.e.\ $\du K\simeq \E^1(K)\cong\E^n(M).$
Hence, $\E^n(M)$ cannot vanish because otherwise $K$ would be
projective.
\end{rem}

\section{APPENDIX: Minimal resolutions}
For lack of a reference we give the proofs of  some basic facts on
minimal resolutions. Let $\La=\zp\kl G\kr $ the completed group
algebra over $\zp$ of a finitely generated pro-$p$-group $G$ and
$k=\La/\M\cong \fp$ its residue class field. We assume that \La\
is Noetherian. For any finitely generated \La-module $M$ we have
the minimal representation
 $$\xymatrix{
   \Lambda^{d_0}\ar@{>>}[r]^{\varphi_0} & M
 }$$
with $d_0=\dim_k M/\M M$ by the Nakayama-Lemma. Proceeding in the
same manner for $\ker(\varphi_0)$ and $d_1=\dim_k
\ker(\varphi_0)/\M\ker(\varphi_0),$ we construct a {\em minimal
free resolution}

$$\xymatrix{
 F_\bullet:   \cdots\ar[r] & {\La^{d_n}_{\phantom{d_n}}}\ar[r]^{\varphi_n} & {\La_{\phantom{d_n}}^{d_{n-1}}}\ar[r]^{\varphi_{n-1}} & {\cdots}\ar[r] & {\La_{\phantom{d_n}}^{d_1}}\ar[r]^{\varphi_1}
   & {\La_{\phantom{d_n}}^{d_0}}\ar[r]^{\varphi_0} & M\ar[r] & 0.}$$

It is easily verified that $F_\bullet$ is determined by $M$ up to
isomorphism of complexes.

\begin{prop}\label{minimal res}
Let $M$ be a finitely generated \La-module and

$$\xymatrix{
 F_\bullet:   \cdots\ar[r] & {F_n}\ar[r]^{\varphi_n} & F_{n-1}\ar[r]^{\varphi_{n-1}} & {\cdots}\ar[r] & F_1\ar[r]^{\varphi_1}
   & F_0\ar[r]  & 0.}$$
a free resolution of $M$. Then the following are equivalent:

\begin{enumerate}
\item $F_\bullet$ is minimal, \item $\varphi_i(F_i)\subseteq\M
F_{i-1}$ for all $i\geq 1,$ \item $\rk_\La(F_i)=\dim_k
\Tor_i^\La(M,k)\fa i\geq 0,$ \item $\rk_\La(F_i)=\dim_k
\Ext^i_\La(M,k)\fa i\geq 0.$
\end{enumerate}
\end{prop}

\begin{proof}
The equivalence of (i) and (ii) follows easily from Nakayama's
lemma. Since $\Tor^\La_i(M,k)=\H_i(F_\bullet\otimes k),$ (iii)
holds if and only if $\varphi_i\otimes k=0$ for all $i\geq 0,$
which is equivalent to (ii). Using
$\Ext_\La^i(M,k)=\H^i(\Hom_\La(F_\bullet,k))$ the equivalence of
(ii) and (iv) follows similarly.
\end{proof}

\begin{cor}\label{mini-res-cor}
Let $M$ be a finitely generated \La-module. Then
 \begin{eqnarray*}
 \pd(M)&=&\max\{i\mid F_i\neq 0\}\\
      &=&\max\{i\mid\Tor^\La_i(M,k)\neq 0\}\\
      &=&\max\{i\mid \Ext^i_\La(M,k)\neq 0\}
\end{eqnarray*}
\end{cor}
\POP

\INPUT{xbib.bib}   
\INPUT{auslander.bbl} 

\bibliographystyle{amsplain}
\bibliography{xbib}

\providecommand{\bysame}{\leavevmode\hbox to3em{\hrulefill}\thinspace}
\providecommand{\MR}{\relax\ifhmode\unskip\space\fi MR }
\providecommand{\MRhref}[2]{%
  \href{http://www.ams.org/mathscinet-getitem?mr=#1}{#2}
}
\providecommand{\href}[2]{#2}
\begin{thebibliography}{10}

\bibitem{aus}
M.~Auslander and M.~Bridger, \emph{Stable module theory}, Memoirs of the AMS,
  vol.~94, AMS, 1969.

\bibitem{bal-how}
P.N. Balister and S.~Howson, \emph{{Note on Nakayama's Lemma for Compact
  $\Lambda$-modules}}, Asian J. Math. \textbf{1} (1997), no.~2, 224--229.

\bibitem{erik1}
J.-E. {Bj\"{o}rk}, \emph{Rings of differential operators}, North-Holland Math.
  Library, vol.~21, North-Holland Publishing Company, 1979.

\bibitem{erik}
\bysame, \emph{{Filtered Noetherian Rings}}, Noetherian rings and their
  applications, Conf. Oberwolfach/FRG 1983, Math. Surv. Monogr., vol.~24, 1987,
  pp.~59--97.

\bibitem{erik2}
\bysame, \emph{{The Auslander condition on Noetherian rings}}, Seminaire
  d'algebre P. Dubreil et M.-P. Malliavin, Proc., Paris/Fr. 1987/88, LNM, vol.
  1404, Springer, 1989, pp.~137--173.

\bibitem{brumer}
A.~Brumer, \emph{Pseudocompact algebras, profinite groups and class
  formations}, J. of Algebra \textbf{4} (1966), 442--470.

\bibitem{bruns}
W.~Bruns and J.~Herzog, \emph{{Cohen-Macaulay rings}}, Cambridge studies in
  advance mathematics, vol.~39, Cambridge University Press, 1993.

\bibitem{cham}
M.~Chamarie, \emph{{Modules sur les anneaux de Krull non commutatfs.}},
  S\'{e}m. d'Alg\`{e}bre P. Dubreil et M.-P. Malliavin 1982, LNM, vol. 1029,
  Springer, 1983, pp.~283--310.

\bibitem{coates}
J.~Coates, \emph{{Fragments of the $\mbox{GL}_2$ Iwasawa theory of elliptic
  curves without complex multiplication.}}, Arithmetic theory of elliptic
  curves. Lectures given at the 3rd session of the Centro Internazionale
  Matematico Estivo (CIME), Cetraro, Italy, July 12-19, 1997., LNM, vol. 1716,
  Springer, 1999, pp.~1--50.

\bibitem{coates-howsonII}
J.~Coates and S.~Howson, \emph{{Euler characteristics and elliptic curves II}},
  J. Math. Soc. Japan \textbf{53} (2001), 175--235.

\bibitem{co-su1}
J.~Coates and R.~Sujatha, \emph{{Structure Theorem}}, private communication
  (2001).

\bibitem{dsms}
J.D. Dixon, M.P.F. du~Sautoy, A.~Mann, and D.~Segal, \emph{Analytic pro-p
  groups}, 1st ed., London Mathematical Society Lecture Note, vol. 157,
  Cambridge University Press, 1991.

\bibitem{fossum}
R.~Fossum, \emph{{Duality over Gorenstein rings}}, Math. Scand. \textbf{26}
  (1970), 177--199.

\bibitem{gelfand-manin}
S.I. Gelfand and Y.I. Manin, \emph{Methods of homological algebra}, Spinger,
  1996.

\bibitem{GoWa}
K.R. Goodearl and R.B. Warfield, \emph{{An Introduction to Noncommutative
  Noetherian Rings}}, LMS Student texts, vol.~16, Cambridge University Press,
  1989.

\bibitem{harris}
M.~Harris, \emph{{p-adic representations arising from descent on Abelian
  varieties}}, Compos. Math. \textbf{39} (1979), 177--245.

\bibitem{hart}
R.~Hartshorne, \emph{{Residues and duality }}, LNM, vol.~20, Springer, 1966.

\bibitem{howson}
S.~Howson, \emph{{Iwasawa theory of Elliptic Curves for $p$-adic Lie
  extensions}}, Ph.D. thesis, University of Cambridge, July 1998.

\bibitem{howson2000}
\bysame, \emph{{Euler Characteristics as Invariants of Iwasawa Modules}},
  preprint (2000).

\bibitem{howson2001}
\bysame, \emph{{Structure of Central Torsion Iwasawa Modules}}, preprint
  (20001).

\bibitem{ja-con}
U.~Jannsen, \emph{{Continuous \'{E}tale Cohomology}}, Math. Ann. \textbf{280}
  (1988), 207--245.

\bibitem{ja-is}
\bysame, \emph{Iwasawa modules up to isomorphism}, Advanced Studies in Pure
  Mathematics \textbf{17} (1989), 171--207.

\bibitem{ja-sps}
\bysame, \emph{{A spectral sequence for Iwasawa adjoints}}, unpublished (1994).

\bibitem{joer}
P.~Joergensen, \emph{{Noncommutative graded homological identities}}, J. Lond.
  Math. Soc., II. Ser. \textbf{57} (1998), no.~2, 336--350.

\bibitem{la}
M.~Lazard, \emph{Groupes analytiques $p$-adiques}, Publ. Math. I.H.E.S.
  \textbf{26} (1965), 389--603.

\bibitem{leva}
T.~Levasseur, \emph{{Grade des modules sur certains anneaux filtres}}, Commun.
  Algebra \textbf{9} (1981), no.~15, 1519--1532.

\bibitem{nsw}
J.~Neukirch, A.~Schmidt, and K.~Wingberg, \emph{Cohomology of number fields},
  Grundlehren der mathematischen Wissenschaften, vol. 323, Springer, 2000.

\bibitem{Ne}
A.~Neumann, \emph{Completed group algebras without zero divisors}, Arch. Math.
  \textbf{51} (1988), 496--499.

\bibitem{ochi-ven}
Y.~Ochi and O.~Venjakob, \emph{{On the structure of Selmer groups over $p$-adic
  Lie extensions}}, to appear in the Journal of Algebraic Geometry.

\bibitem{rot}
J.~J. Rotman, \emph{{An introduction to homological algebra}}, Pure and Applied
  Mathematics, vol.~85, Academic Press, 1979.

\bibitem{schneider2001}
P.~Schneider, \emph{{Notes on modules modulo pseudo-isomorphism}}, private
  communication (2001).

\bibitem{sch-teit}
P.~Schneider and J.~Teitelbaum, \emph{Banach space representations and iwasawa
  theory}, preprint (2000).

\bibitem{swan}
R.~G. Swan, \emph{{Algebraic K-theory}}, LNM, vol.~76, Springer, 1968.

\bibitem{ven-diss}
O.~Venjakob, \emph{{Iwasawa Theory of $p$-adic Lie Extensions}}, Dissertation,
  University of Heidelberg, (2000).

\bibitem{ven}
\bysame, \emph{{Iwasawa Theory of $p$-adic Lie Extensions}}, preprint (2001).

\bibitem{weibel}
C.~A. Weibel, \emph{{An introduction to homological algebra}}, Cambridge
  Studies in Advanced Mathematics, vol.~38, Cambridge University Press, 1995.

\bibitem{wi}
J.~S. Wilson, \emph{Profinite groups}, London Mathematical Society Monographs
  New Series, vol.~19, Oxford University Press, 1998.

\bibitem{zaks}
A.~Zaks, \emph{{Injective dimension of semi-primary rings}}, J. Algebra
  \textbf{13} (1969), 73--86.

\end{thebibliography}
\end{document}